\documentclass{article}

\usepackage{hyperref}  
\usepackage{geometry}  
\usepackage{amssymb, amsbsy}  
\usepackage{mathtools}  
\usepackage{xfrac}  
\usepackage{xspace}  
\usepackage{float}  
\usepackage{subcaption}  
\usepackage{natbib}  
\usepackage[capitalise]{cleveref}  
\usepackage{appendix}  
\usepackage{enumitem}  
\usepackage{multicol}  
\usepackage[noend]{algpseudocode}  
\usepackage{algorithm}  
\usepackage{acro}  
\usepackage{xr}  
\usepackage{bm} 
\usepackage{tabularx}  
\usepackage{etoolbox}  
\usepackage{booktabs}
\usepackage{multirow}
\usepackage{chngcntr}
\usepackage{apptools}
\usepackage{amsfonts}
\usepackage{amssymb}
\usepackage{mathrsfs}
\usepackage{amsmath,amsthm}
\usepackage{authblk}
\usepackage{comment}
\AtAppendix{\counterwithin{lemma}{section}}

\hypersetup{colorlinks=true, allcolors=blue}

\newcommand{\E}{\mathbb{E}}

\newcommand{\mrm}{\mathrm}

\everymath {\displaystyle}

\newcommand{\ruby}[2]{
\leavevmode
\setbox0=\hbox{#1}
\setbox1=\hbox{\tiny #2}
\ifdim\wd0>\wd1 \dimen0=\wd0 \else \dimen0=\wd1 \fi
\hbox{
\kanjiskip=0pt plus 2fil
\xkanjiskip=0pt plus 2fil
\vbox{
\hbox to \dimen0{
\small \hfil#2\hfil}
\nointerlineskip
\hbox to \dimen0{\mathstrut\hfil#1\hfil}}}}
\allowdisplaybreaks[1]

\newcommand{\ms}{\ms}

\newcommand{\new}[3]{\bar{X}_{#3}^{#1, #2}}
\newcommand{\mil}[2]{\bar{X}_{#2}^{#1}}
\newcommand{\EM}[3]{\widetilde{X}_{#3}^{\mrm{EM}, #1, #2}}
\newtheorem{thm}{Theorem}
\newtheorem{lemma}{Lemma}
\newtheorem{prop}{Proposition}

\newtheorem{rem}{Remark}

\newtheorem{cor}{Corollary}

\title{An extended Milstein scheme for 
effective weak approximation of diffusions} 
\author[1]{Yuga Iguchi}
\author[2]{Toshihiro Yamada}
\affil[1]{Department of Statistical Science, University College London, UK.}
\affil[2]{Graduate School of Economics, Hitotsubashi University, Japan.} 

\date{} 

\begin{document}
\maketitle
\begin{abstract} 
We propose a straightforward and effective method for discretizing multi-dimensional diffusion processes as an extension of Milstein scheme. The new scheme is explicitly given and can be simulated using Gaussian variates, requiring the same number of random variables as Euler-Maruyama (EM) scheme. We show that the proposed scheme has a weak convergence rate of one, which is consistent with other classical schemes like EM/Milstein schemes but involves fewer leading-order error terms. Due to the reduction of the error terms, the proposed scheme is expected to provide a more accurate estimation than alternative first-order schemes. We demonstrate that the weak error of the new scheme is effectively reduced compared with EM/Milstein schemes when the diffusion coefficients involve a small parameter. We conduct simulation studies on Asian option pricing in finance to showcase that our proposed scheme significantly outperforms EM/Milstein schemes, while interestingly, we find no differences in the performance between EM and Milstein schemes. 
\end{abstract}

\section{Introduction} 
\label{sec:intro}
This work proposes a simple and effective weak approximation for Stochastic Differential Equations (SDEs) by extending \emph{Milstein scheme}, one of the classical numerical schemes widely exploited in applications. Let $(\Omega, \mathcal{F},  \{\mathcal{F}_t \}_t, \mathbb{P})$ be a filtered probability space and $B_t = (B_t^1, \ldots, B_t^d), \; t \ge 0,$ be the $d$-dimensional standard Brownian motion defined upon the probability space. We consider the $N$-dimensional SDE specified as:
\begin{align}
\begin{aligned} \label{eq:sde}
     d X_t^x = b (X_t^x) dt + \sum_{j = 1}^d \sigma_j (X_t^x) d B_t^j, \quad  X_0^x = x \in \mathbb{R}^N, 
\end{aligned}
\end{align}  
with coefficients $b, \, \sigma_j: \mathbb{R}^N \to \mathbb{R}^N, \, 1 \le j \le d$. 
%
%
Since the solution of SDEs is in general analytically unavailable, some time-discretisation is required to approximately compute the target dynamics. 
The accuracy of the approximation is usually measured by the weak/strong convergence: for a numerical scheme 
$\{\widetilde{X}_{kT/n}^{n, x} \}_{0 \le k \le n}$ starting from a point $x \in \mathbb{R}^N$ with a time $T  > 0$ and a number of discretisation $n \in \mathbb{N}$, 
\begin{gather} 
\Bigl| \mathbb{E} [f (X_T^x)] - \mathbb{E} [f (\widetilde{X}_T^{n, x})]  \Bigr|  = \mathcal{O} (n^{-\alpha});  
\tag{weak error} \\[0.2cm] 
\mathbb{E} \Bigl[ \,  
\max_{0 \le k \le n} 
\bigl| 
X_{kT/n}^x - 
\widetilde{X}_{kT/n}^{n, x} 
\bigr|^2 \, 
  \Bigr]  = \mathcal{O} (n^{-\beta}),  
\tag{strong error} 
\end{gather}
where $f: \mathbb{R}^N \to \mathbb{R}$ is some appropriate test function, and $\alpha, \beta > 0$ corresponds to the order of weak/strong convergence, respectively. For instance, it is well-known that for \emph{Euler-Maruyama (EM) scheme}, $\alpha = 1.0$ and $\beta = 1.0$ (for instance see \cite{kloe:92}). 

As one of the popular discretisation schemes, the {Milstein scheme} has been investigated in many contexts of numerical analysis of SDEs.
It achieves the same rate of weak convergence as the EM scheme ($\alpha = 1.0$) but has an improved rate of strong convergence, i.e., $\beta = 2.0$. The Milstein scheme is defined as: 
for $0 \le k \le n-1$, 
\begin{align} \label{eq:mil}
\begin{aligned}
  \widetilde{X}^{\mrm{Mil}, n , x}_{0} & = x; \\    
  \widetilde{X}^{\mrm{Mil}, n , x}_{(k+1)h}
  & = \widetilde{X}^{\mrm{Mil}, n , x}_{kh}
   + b (\widetilde{X}^{\mrm{Mil}, n , x}_{kh}) h 
   + \sum_{j = 1}^d \sigma_j (\widetilde{X}^{\mrm{Mil}, n , x}_{kh}) \Delta B^{j}_{k}  \\ 
  & \qquad  + \tfrac{1}{2} \sum_{j_1, j_2 = 1}^d 
  g_{j_1 j_2} (\widetilde{X}^{\mrm{Mil}, n , x}_{kh})
  \bigl\{ ( \Delta B^{j_1}_k \Delta B^{j_2}_k - h \times \mathbf{1}_{j_1 = j_2 } ) 
  - A_{j_1j_2, k} 
  \bigr\}, 
\end{aligned} \tag{Milstein}
\end{align}
where we have set: $h = T/ n$, $\Delta B^j_k = B^j_{(k+1)h} - B^j_{kh}$, 
$ \textstyle 
g_{j_1 j_2} (x) 
= \sum_{i = 1}^N \sigma_{j_2}^i (x) \partial_{i} \sigma_{j_1} (x), \, x \in \mathbb{R}^N, 
$
and $A_{j_1j_2, k}$ is the \emph{L\'evy area} defined as: 
\begin{align}
  A_{j_1j_2, k} 
  = \int_{kh}^{(k+1)h} \int_{kh}^{s} 
  d B^{j_1}_u d B^{j_2}_s
  - \int_{kh}^{(k+1)h} \int_{kh}^{s} 
  d B^{j_2}_u d B^{j_1}_s.  
\end{align}
The last term in (\ref{eq:mil}) is obtained from the higher order stochastic Taylor expansion of diffusion coefficients only and contributes to the improvement of the strong convergence, though it cannot be exactly simulated since the law of L\'evy area ($A_{j_1j_2, k}$) is intractable. The Milstein scheme is tractable only if the \emph{commutative condition} holds, i.e., $g_{j_1 j_2}^i (x) - g_{j_2 j_1}^i (x) = 0, \, 1 \le i \le N, \, 1 \le j_1, j_2 \le d$ for any $x \in \mathbb{R}^N$. Precisely, under the commutative condition, the terms involving $A_{j_1 j_2, k}$ are cancelled out due to its anti-symmetric property. However, many important SDE models in applications, e.g., \emph{stochastic volatility models} in financial engineering, do not satisfy the commutative condition. Thus, in practice, use of \emph{truncated Milstein scheme}, i.e., the scheme (\ref{eq:mil}) with $A_{j_1 j_2, k}$ being replaced with $0$, will be a one of realistic options. The scheme is defined as: 
\begin{align} \label{eq:T-mil}
\begin{aligned}
\widetilde{X}^{\mrm{TMil}, n , x}_{0} & = x;   \\ 
\widetilde{X}^{\mrm{T}\mrm{Mil}, n , x}_{(k+1)h}
& = \widetilde{X}^{\mrm{TMil}, n , x}_{kh}
+ b (\widetilde{X}^{\mrm{TMil}, n , x}_{kh}) h 
+ \sum_{j = 1}^d \sigma_j (\widetilde{X}^{\mrm{TMil}, n , x}_{kh}) \Delta B^{j}_{k}  \\ 
& \qquad  + \tfrac{1}{2} \sum_{j_1, j_2 = 1}^d 
g_{j_1 j_2} (\widetilde{X}^{\mrm{TMil}, n , x}_{kh})
( \Delta B^{j_1}_k \Delta B^{j_2}_k - h \times \mathbf{1}_{j_1 = j_2 }). 
\end{aligned} \tag{T-Milstein}
\end{align} 
The truncated scheme achieves the first order weak convergence, $\alpha = 1.0$. However, the removal of $A_{j_1 j_2, k} $ leads to deterioration of the strong convergence rate under the non-commutative condition and thus $\beta = 1.0$. \cite{yamada23} showed that for elliptic diffusions containing small parameters in both drift and diffusion coefficients, the truncated Milstein scheme can effectively reduce the weak error compared with the EM scheme. 
The truncated Milstein scheme is also exploited in the framework of multi-level Monte-Carlo (MLMC).  For instance, \cite{giles14} developed an antithetic multi-level Monte-Carlo based upon the truncated Milstein scheme.  

One of the main focuses in this work is to investigate the (truncated) Milstein scheme from the perspective of weak approximation, motivated by the following question: 
\textit{Can the Milstein scheme or the truncated Milstein scheme always produce a more accurate weak approximation than the EM scheme with the aid of higher order stochastic Taylor expansion of the diffusion coefficients?} 
The analytic result and numerical experiments provided later in this paper demonstrate that using the Milstein scheme does not necessarily have some advantages over the Euler-Maruyama scheme in weak approximation. In particular, this is shown by simulation studies in Section \ref{sec:numerics}, where we find no significant differences in the performance of those two schemes. 
Then, instead of using these first order discretisation schemes, one might consider higher order weak approximation such as KLNV method \citep{gyur10, gyur11, lit12, NV08}, Malliavin weight approach \citep{OS21, igu21, ny19, ty16, y19, yy20},  or weak second order sampling schemes \citep{mil21, talay84}. However, to improve the weak convergence, those schemes require extra computational efforts to take the L\'evy area into account. For instance, most of the above methods require simulation of $D$-i.i.d. random variables for the local (one-step) transition, where $D$ is larger than $d$, i.e., the dimension of driving Brownian motion. Thus, when $n$-times discretisation is introduced, the total number of generated random variables becomes $nD$, larger than $nd$ required for EM and (truncated) Milstein schemes. 

The central objective of this work is to propose a simple and effective first order weak approximation for a wide class of diffusions to outperform the Euler-Maruyama and the (truncated) Milstein scheme with a similar computational cost. In brief, the scheme we propose in this article features the following:
\begin{itemize}[leftmargin=0.4cm]
  \item It is always \emph{explicit} even if the commutative condition does not hold.    
  \item It extends scheme (\ref{eq:T-mil}) and involves additional terms from stochastic Taylor expansion of the drift coefficients. 
  \item The total number of random variables to simulate the scheme is $nd$ for $n$-times discretization, which is the same as EM and (truncated) Milstein schemes. 
\end{itemize} 

We refer to our proposed scheme as \emph{extended Milstein scheme}. We will show (under some conditions on the SDE's coefficients) that the extended Milstein scheme achieves a first order weak approximation when the test function is bounded and measurable (non-smooth). 
In particular, we derive the analytic weak error expansion of the new scheme and compare it with that of EM scheme and (truncated) Milstein scheme. Then, for a class of hypo-elliptic diffusions with diffusion coefficients containing a small parameter $\varepsilon \in (0,1)$, we quantitatively show that the discretization bias by the extended Milstein scheme is smaller than that of other first order schemes with the help of the small parameter $\varepsilon$. 
Finally, we present simulation studies showcasing that the extended Milstein scheme indeed significantly reduces the weak approximation bias compared with the other classical numerical schemes.  
\\

\noindent
\textbf{Notation.}
%
%
We denote by $\mathscr{B}_b (\mathbb{R}^m), \, m \ge 1$, the space of bounded and measurable functions $f : \mathbb{R}^m \to \mathbb{R}$. We write $C_b^\infty (\mathbb{R}^{m}; \mathbb{R}^{n}), \, m, n \ge 1$ as the space of smooth functions $\varphi : \mathbb{R}^{m} \to \mathbb{R}^{n}$ with bounded derivatives of any order. 
We set $a = \sigma \sigma^\top$ and write the generator of the SDE (\ref{eq:sde}) as: 
\begin{align} \label{eq:L}
     \mathcal{L} \equiv \sum_{i = 1}^N b^i (\cdot) \, \partial_i + \frac{1}{2} \sum_{i_1, i_2 = 1}^N 
     a^{i_1 i_2} (\cdot) \,  \partial_{{i_1} {i_2}}
 \end{align}
Throughout the paper, let $T >0$, and we frequently use the notation of the test function $f \in \mathscr{B}_b (\mathbb{R}^N)$. We define the uniform norm $\| \cdot \|_\infty$ as: $\textstyle \| f \|_\infty = \sup_{x \in \mathbb{R}^N} | f (x) |$ for $f \in \mathscr{B}_b (\mathbb{R}^N)$. We write $u (s,x) = \E [ f (X_{T-s}^x) ], \, (s, x) \in [0,T] \times \mathbb{R}^N$ as the solution to the parabolic partial differential equation (PDE) specified as: 
\begin{align}
\begin{aligned} 
   \partial_t u (t,x) + \mathcal{L} u (t,x) & = 0,
  \ \ (t,x) \in [0,T)  \times \mathbb{R}^N, \\ 
   u (T,x)  & = f (x).  
\end{aligned}       
\end{align} 
%
We define a function $\widetilde{b} : \mathbb{R}^N \to \mathbb{R}^N$ as: 
\begin{align} \label{eq:drift_st}
\widetilde{b}^i (x) = b^i (x) - \tfrac{1}{2} \sum_{k = 1}^d 
\sum_{j= 1}^N \sigma_k^j (x) \partial_j \sigma_k^i (x), \quad x \in \mathbb{R}^N, \quad  1 \le i \le N, 
\end{align}
which is the drift coefficient of the SDE when the It\^o-type SDE (\ref{eq:sde}) is written in the Stratonovich form. We identify the coefficients $\widetilde{b}, \sigma_j, \, 1 \le j \le  d$ with the following vector fields: 
\begin{align} \label{eq:vf}
    \widetilde{L}_0 & = \sum_{i = 1}^N \widetilde{b}^i (\cdot) \, 
    \partial_i, \quad    
    L_j  =  \sum_{i = 1}^N \sigma_j^i (\cdot) \,  \partial_i, \quad  
    1 \le j  \le d.  
\end{align}
For two vector fields $V, W$, the Lie bracket is defined as: 
\begin{align} \label{eq:lb}
    [V, W] = VW - WV. 
\end{align}
\section{Preliminaries} \label{sec:pre}
We introduce two main conditions for the SDE (\ref{eq:sde}), and then briefly review the work of \citep{bally96}, which provided an analytic result for weak approximation by Euler-Maruyama (EM) scheme under the above conditions. Based upon the analytic error expansion of EM scheme in the weak sense, we will compare the discretisation biases induced by the EM scheme and our new scheme proposed later in the next section. 
\subsection{Conditions for the SDE} 
We define sets of vector fields constructed from the Lie bracket (\ref{eq:lb}) as:
\begin{gather*}
\Sigma_0 = \bigl\{ L_1, \ldots, L_d \bigr\}, 
\quad \Sigma_m = \Bigl\{ 
  \bigl\{ [\widetilde{L}_0, V], \, [L_j, V] \bigr\} 
  \,  :  V \in \Sigma_{m-1}, \, 1 \le j \le d \Bigr\}, \quad  m \ge 1. 
\end{gather*} 
For the SDE (\ref{eq:sde}), or the Stratonovich-SDE with the drift (\ref{eq:drift_st}), we assume the following two conditions under which \cite{bally96} studied the weak error of Euler-Maruyama scheme when the test function is assumed to be bounded and measurable:
\begin{enumerate}
\item \label{con:coeff}
$b, \, \sigma_j, \, j = 1, \ldots, d,$ are infinitely differentiable, and their derivatives of any order are bounded.
\item \label{con:hor}
Uniform H\"ormander's condition holds: there exists an integer $M \ge 0$ such that 
\begin{align}
\inf_{x \in \mathbb{R}^N} 
\inf_{\substack{\xi \in \mathbb{R}^N \\
\mrm{s.t.} \|\xi\| = 1}} 
\sum_{0 \le i \le M} \sum_{V \in \Sigma_i }
\langle V(x), \xi \rangle^2 > 0. 
\end{align}  
\end{enumerate}
We note that under the condition (\ref{con:hor}), the law of $X_t^x$ is absolutely continuous with respect to (w.r.t.) the Lebesgue measure for any $x \in \mathbb{R}^N$ and $t > 0$. When (\ref{con:hor}) holds with $M =0$, the condition is typically interpreted as the uniformly elliptic condition, and thus the matrix $a := \sigma \sigma^\top$ is positive definite uniformly in the state $x \in \mathbb{R}^N$. Furthermore, together with the condition (\ref{con:coeff}), the condition (\ref{con:hor}) leads to the existence of a smooth Lebesgue density of the law of $X_t^x$, see, e.g. \cite{nua06}.
\subsection{Weak approximation by the Euler-Maruyama scheme -- Review}
We review the theoretical result for weak approximation by the EM scheme for hypo-elliptic diffusions. The EM scheme is defined as: for $x \in \mathbb{R}^N$, $n \ge 1$ and $h := T/n$, 
\begin{align}
\begin{aligned} \label{eq:em}
\EM{n}{x}{0} & = x; \\ 
\EM{n}{x}{(k+1)h} & = \EM{n}{x}{kh} 
+ b \bigl( \EM{n}{x}{kh}  \bigr) h 
+ \sum_{j = 1}^d \sigma (\EM{n}{x}{kh}) \,\Delta B_k^j, \quad 0 \le k \le n-1, 
\end{aligned}  
\tag{EM}
\end{align}
Under conditions (\ref{con:coeff})--(\ref{con:hor}), \cite{bally96} showed that the EM scheme achieves the first order weak convergence:
\begin{thm}[\cite{bally96}] \label{thm:wa_em}
Let $T > 0$, $x \in \mathbb{R}^N$ and 
$f \in \mathscr{B}_b (\mathbb{R}^N)$. Under the conditions (\ref{con:coeff})--(\ref{con:hor}), it holds that: 
\begin{align} \label{eq:wa_em}
    \E [f (X_T^x)] - \E [f (\EM{n}{x}{T})] 
    = \frac{T}{n} \,  C_f^{\mrm{EM}} (T, x)   + \frac{1}{n^2} \, R_f^{\mrm{EM}} (T,x),  
\end{align}
with the remainder terms $C_f^{\mrm{EM}} (T, x)$ and $R_f^{\mrm{EM}} (T, x)$ specified as follows: the leading order term is given by $\tfrac{T}{n} C_f^{\mrm{EM}} (T, x)$ with 
\begin{align*} 
C_f^{\mrm{EM}} (T, x) = \int_0^T  \mathbb{E} \bigl[ \Phi^{\mrm{EM}}  (s, X_s^x) \bigr] ds, 
\end{align*}
where $\Phi^{\mathrm{EM}} : [0,\infty) \times \mathbb{R}^N \to \mathbb{R}$ is defined as:   
\begin{align}
\begin{aligned} \label{eq:err_em_original}
& \Phi^{\mathrm{EM}} (t,x) 
=  - \tfrac{1}{2} \sum_{i, j = 1}^N b^i (x) b^j (x) \partial_{i j} u (t, x) 
-  \tfrac{1}{2} \sum_{i, j, k  = 1}^N b^i (x)  a^{jk}(x) \partial_{ijk} u (t,x) \\ 
& \quad - \tfrac{1}{8} \sum_{i,j,k,l = 1}^N a^{ij}(x)  a^{kl}(x) \partial_{ijkl} u (t, x) 
{ -  \tfrac{1}{2} {\partial_t^2} u (t, x)} 
{- \sum_{i = 1}^N  b^i (x) {\partial_t} \partial_i u(t,x) 
- \tfrac{1}{2} \sum_{i, j = 1}^N a^{ij} (x)  {\partial_t} \partial_{ij} u(t,x)}. 
\end{aligned} 
\end{align}
Also, there exist real numbers $q, Q > 0$ and a non-decreasing function $K(\cdot)$ such that
\begin{align} 
 | C_f^{\mrm{EM}} (T, x) | + | R_f^{\mrm{EM}} (T, x) | 
 \leq K(T) \| f \|_\infty \frac{1 + |x|^Q}{T^q}.  
\end{align}
\end{thm}
%
%
%
Our objective is to compare the weak approximation error by EM scheme with that by the extended Milstein scheme proposed later in (\ref{eq:ext-Mil}). To this end, we here provide an expression for $\Phi^{\mrm{EM}}$ without requiring the partial derivatives of $u$ w.r.t. the time variable. For simplicity of the notation, we write: 
\begin{align}
L_0 \equiv \mathcal{L} = \sum_{i = 1}^N b^i (\cdot) \partial_i + \tfrac{1}{2} \sum_{1 \le i,j \le N} a^{ij} (\cdot) \partial_{ij}.  
\end{align}
Making use of $\partial_t u = - \mathcal{L} u$, we obtain the following result whose proof is postponed to Appendix \ref{appendix:Phi_em}. 
\begin{lemma} \label{lemm:Phi_em}
The function $\Phi^{\mathrm{EM}} : [0,\infty) \times \mathbb{R}^N \to \mathbb{R}$ is given as: 
\begin{align} \label{eq:Phi_em}
\Phi^{\mathrm{EM}} (s ,x) = \Phi_1 (s ,x) + \Phi_2  (s ,x) + \Phi_3 (s ,x), 
\quad (s, x) \in [0,\infty) \times \mathbb{R}^N, 
\end{align}
with  
\begin{align}
 \Phi_{1}  (s, x)  
& \equiv \tfrac{1}{2} \sum_{i = 1}^N L_0 b^i (x) \partial_i u (s, x) 
+ \tfrac{1}{2} \sum_{i,j = 1}^N \sum_{m = 1}^d 
\sigma_{m}^i (x) \bigl\{  L_{m} b^j (x) + L_0 \sigma_m^j (x) \bigr\} \partial_{ij} u (s,x),  
\label{eq:Phi_1} \\[0.2cm] 
\Phi_{2}  (s, x) 
& \equiv \tfrac{1}{2} \sum_{i, j, k}^N \sum_{m_1, m_2 = 1}^d 
L_{m_1} \sigma_{m_2}^i (x) \sigma_{m_1}^j (x) \sigma_{m_2}^k (x) \partial_{ijk} u (s, x)  \nonumber \\ 
& \qquad  + \tfrac{1}{8} \sum_{i,j = 1}^N \sum_{m_1, m_2 =1}^d 
L_{m_1} \sigma_{m_2}^i (x) \bigl\{ L_{m_1} \sigma_{m_2}^j (x) + L_{m_2} \sigma_{m_1}^j (x) \bigr\} \partial_{ij} u (s, x), 
\label{eq:Phi_2} \\
\Phi_3 (s, x) 
& \equiv 
\tfrac{1}{8} \sum_{i, j = 1}^N \sum_{m_1 , m_2 = 1}^d 
L_{m_1} \sigma_{m_2}^i (x) \, [L_{m_1}, L_{m_2}]^j (x) 
\, \partial_{ij} u (s,x). 
\label{eq:Phi_3}
\end{align}
\end{lemma}
Thus, due to Theorem \ref{thm:wa_em} and Lemma \ref{lemm:Phi_em}, the weak approximation error invoked by the EM scheme involves leading order term of size $\mathcal{O}(T/n)$ that is specified as: 
\begin{align} \label{eq:err_em}
 \frac{T}{n} \times  C_f^{\mrm{EM}} (T,x) 
 = \frac{T}{n}  \times \int_0^T 
 \E \bigl[ \Phi_1 (s, x) + \Phi_2 (s, x) 
 + \Phi_3 (s, x)  \bigr] ds. 
\end{align} 
\section{Extended Milstein scheme and weak error analysis}
In this section, we propose the extended Milstein scheme mentioned in Section \ref{sec:intro} so that it produces a more effective weak approximation than other classical numerical schemes. We then show the analytic weak error expansions for the EM, the (truncated) Milstein and the extended Milstein scheme in the manner of Theorem \ref{thm:wa_em} as the core result in this paper. In the last subsection, we apply the main result for a class of hypo-elliptic diffusions with a small parameter in the diffusion coefficients to clarify an advantage of the proposed scheme over other classical schemes.   
\subsection{Extended Milstein scheme}
We introduce a new discretisation scheme as an extension of the truncated Milstein scheme. For simplicity of notation, we write $b \equiv \sigma_0$. Let $T > 0$, $x \in \mathbb{R}^N$ and $h = T/n$ with $n \in \mathbb{N}$. Then, we propose the extended Milstein scheme as follows: for $0 \le k \le n-1$,  
\begin{align} \label{eq:ext-Mil}
\begin{aligned}
& \new{n}{x}{0} = x; \\ 
& \new{n}{x}{(k+1)h}  = \new{n}{x}{kh}  +  
\sum_{j = 0}^d \sigma_j (\new{n}{x}{kh}) \Delta B_k^j 
+ \sum_{0 \le j_1, j_2 \le d} L_{j_1} \sigma_{j_2} (\new{n}{x}{kh}) \times \frac{1}{2}
\Bigl\{ \Delta B^{j_1}_{k} \Delta B^{j_2}_{k} 
- h \times \mathbf{1}_{j_1 = j_2 \neq 0}
\Bigr\},
\end{aligned} 
\end{align}
where we interpret $\Delta B_k^0 = h$.
We notice that the scheme (\ref{eq:ext-Mil}) is explicitly given and simulated by $d$-dimensional Brownian increments only for each iteration. Also, it includes terms from the stochastic Taylor expansion of the drift function $\sigma_0$ and $L_0 \sigma_{j_2}$. 

We have the following result on the proposed scheme (\ref{eq:ext-Mil}): 
\begin{thm} \label{thm:wa_mil} 
Let $T > 0$, $x \in \mathbb{R}^N$ and $f \in \mathscr{B}_b (\mathbb{R}^N)$. Under the conditions (\ref{con:coeff})--(\ref{con:hor}), it holds that:
\begin{align} \label{eq:new_wa}
    \E [f (X_T^x)] - \E [ f (\new{n}{x}{T})] 
    = \frac{T}{n} \, C_f (T, x)  
    + \frac{1}{n^2} \,  R_f (T,x),  
\end{align}
with the terms $C_f (T, x)$ and $R_f (T, x)$ specified as follows: the first term of the right hand side of (\ref{eq:new_wa}) is the leading order term with 
\begin{align} \label{eq:err_new}
C_f (T, x) \equiv 
\int_0^T  \mathbb{E} \bigl[ \Phi_3  (s, X_s^x) \bigr] ds, 
\end{align}
where $\Phi_3$ is given in (\ref{eq:Phi_3}). 
Also, there exist real numbers $q, Q > 0$ and a non-decreasing function $K(\cdot)$ such that
\begin{align}  \label{eq:bd_weak_error}
| C_f (T, x) | + | R_f (T, x) |
\leq K(T) \| f \|_\infty \frac{1 + |x|^Q}{T^q}.  
\end{align}
 \end{thm}
The proof of Theorem \ref{thm:wa_mil} is contained in Section \ref{sec:pf_main}. Notice that the leading error term by the extended Milstein scheme (\ref{eq:ext-Mil}) is given as:
\begin{align} \label{eq:new_leading}
\frac{T}{n} C_f (T,x) = \frac{T}{n} \int_0^T \E [\Phi_3 (s, x)] ds.
\end{align}   
%
Upon consideration of Theorem \ref{thm:wa_em}, both the Euler-Maruyama scheme and the extended Milstein scheme share the same rate of weak convergence $\mathcal{O}(n^{-1})$, but the extended Milstein scheme (\ref{eq:ext-Mil}) invoke the fewer error terms than Euler-Maruyama scheme does: See $\textstyle{\tfrac{T}{n} \, C_f^{\mrm{EM}} (T,x)}$ and $\textstyle{\tfrac{T}{n} \, C_f (T,x)}$  given in (\ref{eq:err_em}) and (\ref{eq:new_leading}), respectively.  
\begin{rem} \label{rem:err}
It is shown that the leading weak error terms of the (intractable) Milstein scheme and the (tractable) truncated Milstein scheme are given as: 
\begin{align}
\tfrac{T}{n} \int_0^T  \E \bigl[ \Phi_1 (s, X_s^x) \bigr] ds,  \qquad 
\tfrac{T}{n}  \int_0^T  \E \bigl[ 
  \Phi_1 (s, X_s^x)
  + \Phi_3 (s, X_s^x)  \bigr] ds,   
\end{align}
respectively. The above error terms and Theorem \ref{thm:wa_mil} imply that the appearance of the term $\Phi_3$ results from removing the L\'evy area in the definition of the truncated Milstein/extended Milstein scheme. 
Furthermore, in the case of the proposed scheme (\ref{eq:ext-Mil}), the weak error term $\textstyle{\frac{T}{n} \int_0^T \E [\Phi_1 (s, x)] ds}$ does not appear. This is because the scheme (\ref{eq:ext-Mil}) involves the following terms in its definition: 
\begin{align} 
L_0 \sigma_0  (\new{n}{x}{kh}) \times \frac{h^2}{2}, 
\quad 
\sum_{1 \le j \le d} \Bigl\{ L_0 \sigma_j  (\new{n}{x}{kh})  
+ L_j \sigma_0  (\new{n}{x}{kh}) \Bigr\} 
\times \frac{1}{2} \Delta B_k^j h, \quad 0 \le k \le n-1. 
\end{align}
\end{rem} 
\begin{table}
  \caption{{Comparison of first order schemes}} 
  \centering 
   \label{table:comparison} 
  \begin{tabular}{cccc} 
  \toprule 
  \\[-10pt] 
  Scheme  & Leading error term & Requirement of L\'evy area
  \\
  \midrule  
  \ref{eq:em}
  & 
  $\tfrac{T}{n} \sum_{1 \le i \le 3} \int_0^T \mathbb{E} [\Phi_i (s, X_s^x)] ds$ 
  & No 
  \\[0.1cm]
  \ref{eq:mil}
  & 
  $\tfrac{T}{n} \int_0^T \mathbb{E} [\Phi_1 (s, X_s^x)] ds$ 
  & Yes 
  \\[0.1cm]
  \ref{eq:T-mil}
  & 
  $\tfrac{T}{n} 
  \int_0^T \mathbb{E} [\Phi_1 (s, X_s^x)
  + \Phi_3 (s, X_s^x)] ds$ 
  & No 
  \\[0.1cm]
  New (\ref{eq:ext-Mil})
  & 
  $\tfrac{T}{n} 
  \int_0^T \mathbb{E} [\Phi_3 (s, X_s^x)] ds$ 
  & No
  \\ 
  \bottomrule   
  \end{tabular}
  \end{table} 
\subsection{Hypo-elliptic diffusions with a small diffusion parameter}
\label{sec:small_diff}
Based on the main result (Theorem \ref{thm:wa_mil}),  we analytically demonstrate the advantage of the proposed scheme (\ref{eq:ext-Mil}) over the Euler-Maruyama/Milstein schemes in terms of weak approximation. In particular, when the diffusion coefficients involve a small parameter $\varepsilon \in (0,1)$, which often appears in applications, we quantitatively show that the leading weak error term of the extended scheme (\ref{eq:ext-Mil}) is smaller than that of EM and (truncated) Milstein scheme by incorporating $\varepsilon$ into the error bound. 
To observe this, we introduce the following hypo-elliptic diffusion:  
\begin{align} \label{eq:small_diffusion}
\begin{aligned} 
d X_t^{x, \varepsilon}  = 
\begin{bmatrix}
d X_{R, t}^{x, \varepsilon} \\[0.1cm] 
d X_{S, t}^{x, \varepsilon} 
\end{bmatrix}
= 
\begin{bmatrix}
b_{R} (X_t^{x, \varepsilon})  \\[0.1cm] 
b_{S} (X_t^{x, \varepsilon})
\end{bmatrix} dt 
+ \sum_{j = 1}
^d\begin{bmatrix}
\epsilon \, \sigma_j (X_t^{x, \varepsilon})  \\[0.1cm]
\mathbf{0}_{N_S}
\end{bmatrix} dB_t^j,  \quad X_0^{x, \varepsilon} = (X_{R, 0}^{x, \varepsilon}, X_{S, 0}^{x, \varepsilon}) = (x_R, x_s) \in \mathbb{R}^N, 
\end{aligned}
\end{align} 
where $\varepsilon \in (0,1)$ and $X_{R,t}^{x, \varepsilon} \in \mathbb{R}^{N_R}, X_{S,t}^{x, \varepsilon} \in \mathbb{R}^{N_S}$ with integers $N_R, N_S \ge 1$ such that $N_R + N_S = N$. In the above, the coefficients are specified as: 
\begin{gather*}
b_R : \mathbb{R}^N \to \mathbb{R}^{N_R}, 
\quad b_S : \mathbb{R}^N \to \mathbb{R}^{N_S}, 
\quad \sigma_j : \mathbb{R}^N \to \mathbb{R}^{N_R}, \quad 1 \le j \le d.  
\end{gather*}
When considering the model (\ref{eq:small_diffusion}), we write $u^\varepsilon (s, x), \, (s, x) \in (0, \infty) \times \mathbb{R}^d$ instead of $u (s, x)$ to emphasise the dependence of the parameter $\varepsilon$. We then introduce the following result whose proof is provided in Appendix \ref{appendix:small_diff_bd}. 
\begin{prop} \label{prop:small_diff_bd}
Let $T > 0$, $\alpha \in \{1, \ldots, N \}^k, k \in \mathbb{N}$ and $g \in \mathbb{R}^N \to \mathbb{R}$ be a smooth function with polynomial growth. 
Work under the conditions (\ref{con:coeff}) and (\ref{con:hor}) with $M = 1$.  Then,  there exist a non-decreasing function $K(\cdot)$ and constants  $q, Q > 0$ independent of $ \varepsilon \in (0,1), n, N$ such that
\begin{align}
\Bigl| \int_0^T \E \bigl[ g (X_s^{x, \varepsilon}) 
\partial_{\alpha} u^\varepsilon (s, X_s^{x, \varepsilon}) \bigr]  ds  \Bigr| \le 
K(T)  
\frac{\|  f \|_\infty}{T^q \varepsilon^k} \times (1 + |x|^Q).  
\end{align}
\end{prop}
We also adjust the notation of the EM/Milstein/truncated Milstein/extended Milstein scheme applied to the SDE (\ref{eq:small_diffusion}) as $\widetilde{X}_{kh}^{\mrm{EM}, n, x, \varepsilon} / \widetilde{X}_{kh}^{\mrm{Mil}, n, x, \varepsilon} /  \widetilde{X}_{kh}^{\mrm{TMil}, n, x, \varepsilon} / \bar{X}_{kh}^{n, x, \varepsilon}$, respectively. 
Application of Proposition \ref{prop:small_diff_bd} to the leading error terms for the first order schemes presented in Table \ref{table:comparison}, 
together with Theorems \ref{thm:wa_em} and \ref{thm:wa_mil}, leads to the following result.  
\begin{cor} \label{cor:sharp_bd}
Let $ T > 0$ , $x \in \mathbb{R}^N$, $\varepsilon \in (0,1)$ and $n \in \mathbb{N}$.  Assume conditions (\ref{con:coeff}) and (\ref{con:hor}) with $M=1$ hold. 
\begin{enumerate}[leftmargin=0.4cm]
\item[1.] Let $f_{\mrm{Lip}} : \mathbb{R}^N \to \mathbb{R}$ be Lipschitz continuous. There exist positive constants $C_{1}^{w}, q_1^{w}, \, w \in \{\mrm{EM}, \mrm{T Mil}, \mrm{Mil}, \mrm{New} \}$,  independent of $n, \varepsilon $ and $x$ such that
\begin{align*} 
 \Bigl| \E \bigl[ f_{\mrm{Lip}} (X_T^{x, \varepsilon})  \bigr] -
\E \bigl[ f_{\mrm{Lip}} ( \widetilde{X}^{w, n, x, \varepsilon}_T ) \bigr]  
 \Bigr| & \le C_{1}^{w} \times \frac{\|\nabla f_{\mrm{Lip}} \|_\infty}{n} (1 + |x|^{q_1^{w}}), \qquad w \in \{\mrm{EM}, \mrm{TMil}, \mrm{Mil} \}; \\[0.2cm] 
\Bigl| \E \bigl[ f_{\mrm{Lip}} (X_T^{x, \varepsilon})  \bigr] - 
\E \bigl[ f_{\mrm{Lip}} ( \new{n}{x, \varepsilon}{T} ) \bigr]   \Bigr| 
& \le C_{1}^{\mrm{New}} \varepsilon^3 \times \frac{\|\nabla f_{\mrm{Lip}} \|_\infty}{n} (1 + |x|^{q_1^{\mrm{New}}}).   
\end{align*}
\item[2.] Let $f \in \mathscr{B}_b (\mathbb{R}^N)$.  There exist positive constants $C_{2}^{w}, q_2^{w}, \, w \in \{\mrm{EM}, \mrm{T Mil}, \mrm{Mil}, \mrm{New} \}$, independent of $n, \varepsilon $ and $x$ such that
\begin{align*} 
\Bigl| \E \bigl[ f  (X_T^{x, \varepsilon})  \bigr] -
\E \bigl[ f ( \widetilde{X}^{w, n, x, \varepsilon}_T ) \bigr] \Bigr|
& \le C_{2}^{w} \times  \frac{\| f \|_\infty}{n \varepsilon} (1 + |x|^{q_2^{w}}), 
\quad w \in \{\mrm{EM}, \mrm{TMil}, \mrm{Mil} \}; 
\\[0.2cm]  
\Bigl| \E \bigl[ f (X_T^{x, \varepsilon})  \bigr] -
\E \bigl[ f  ( \new{n}{x, \varepsilon}{T} ) \bigr]   \Bigr| 
& \le C_{2}^{\mrm{New}} \varepsilon^2 \times \frac{\| f \|_\infty}{n} (1 + |x|^{q_2^{\mrm{New}}}).   
\end{align*} 
\end{enumerate}
\end{cor}

Corollary \ref{cor:sharp_bd} demonstrates that the extended scheme (\ref{eq:ext-Mil}) has a better upper bound with the help of $\varepsilon \in (0,1)$, compared with other schemes. The differences in the upper bounds result from that EM, truncated Milstein and Milstein scheme involve the term $\textstyle \int_0^T \mathbb{E} \bigl[ \Phi_1 (s, X_s^{x, \varepsilon}) \bigr] ds$ in the leading error, while the term is eliminated in the case of the extended Milstein scheme (see also Remark \ref{rem:err}). 
Furthermore, due to the presence of $\varepsilon^\alpha, \, \alpha > 0$ in the bound, the proposed scheme can behave nearly as a second order weak approximation up to some number of discretisation $n \in \mathbb{N}$ given a sufficiently small $\varepsilon \in (0,1)$. 
Thus the extended Milstein scheme provides an asymptotic approximation which is slightly different to the discretisation schemes or expansions in \citep{igu21, ty12, ty16}. 
This is also observed in the simulation studies in the following section. 
\begin{rem}
  The interpretation of Corollary \ref{cor:sharp_bd} is as follows. When the SDE is driven by a small noise, the model is mainly dominated by the drift function and thus, it is critical to include the higher order expansion terms of the drift function to increase the precision of approximation. 
  We iterate here that the three schemes (EM/truncated Milstein/Milstein) do not exploit such higher order terms from the drift function while the extended Milstein scheme does.    
\end{rem}
\section{Simulation studies} 
\label{sec:numerics}
\subsection{Asian call option pricing with Black-Scholes model}
\begin{figure}[h]
  \centering
  \begin{subfigure}[b]{0.62 \textwidth}
      \centering
      \includegraphics[width=\textwidth]{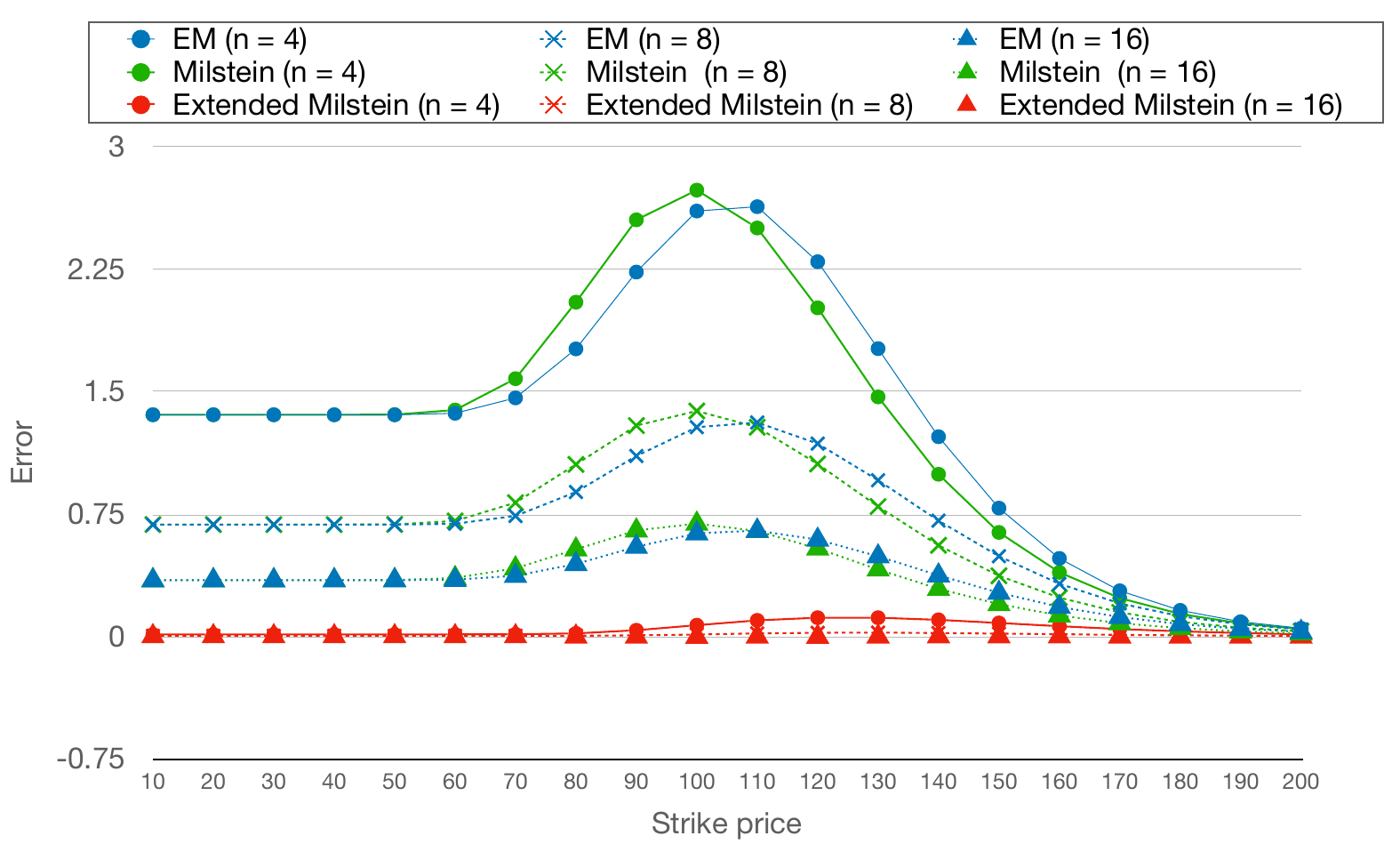}
     \caption{Error of estimation ($\sigma = 0.4$).} %
      \label{fig:asian_bs_err_1}
  \end{subfigure}
  \hspace{0.5cm}
  \begin{subfigure}[b]{0.32 \textwidth}
      \centering
      \includegraphics[width=\textwidth]{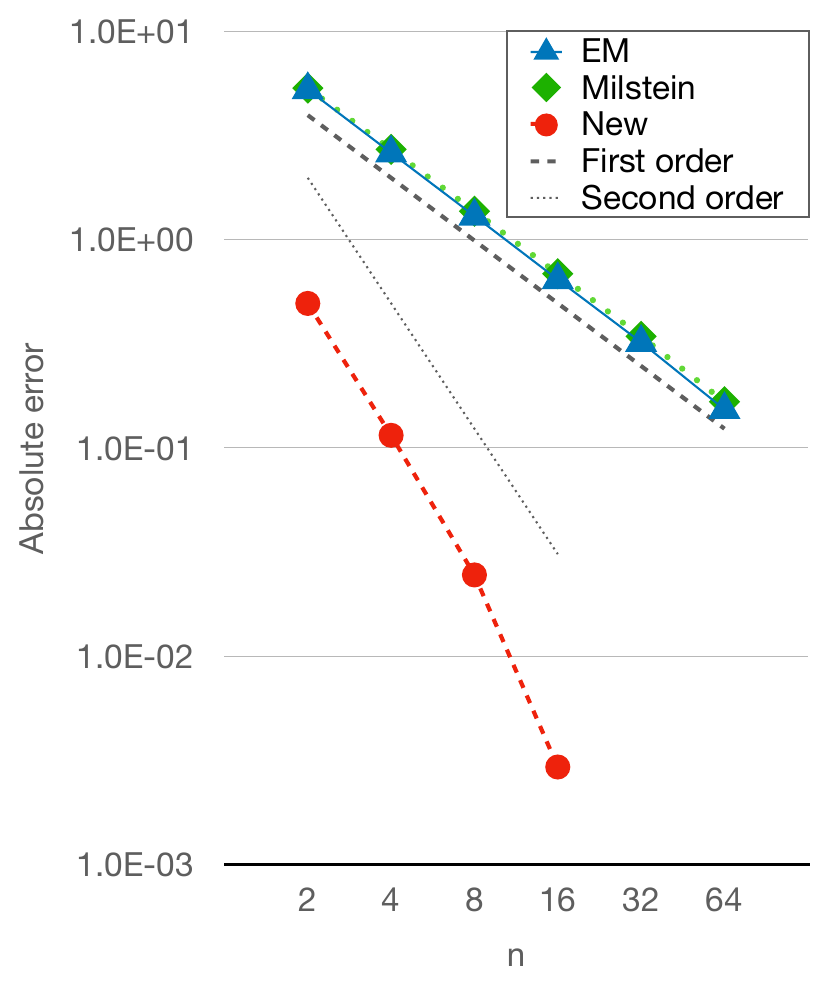}
     \caption{Weak convergence ($\sigma = 0.4$).} 
      \label{fig:asian_bs_conv_1}
  \end{subfigure} 
  \\[1cm] 
  \begin{subfigure}[b]{0.62 \textwidth}
    \centering
    \includegraphics[width=\textwidth]{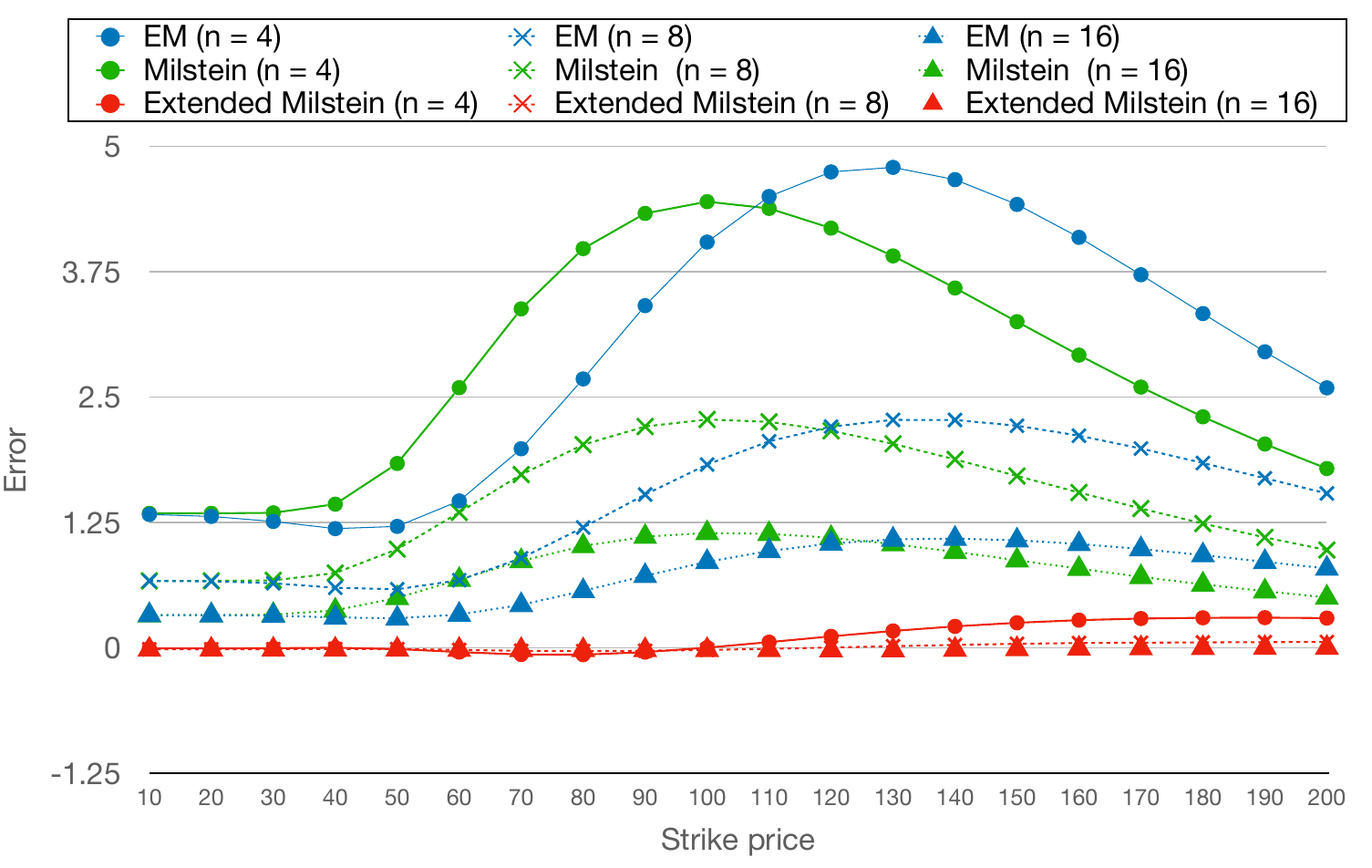}
   \caption{Error of estimation ($\sigma = 0.8$).} %
    \label{fig:asian_bs_err_2}
\end{subfigure}
\hspace{0.5cm}
\begin{subfigure}[b]{0.33 \textwidth}
    \centering
    \includegraphics[width=\textwidth]{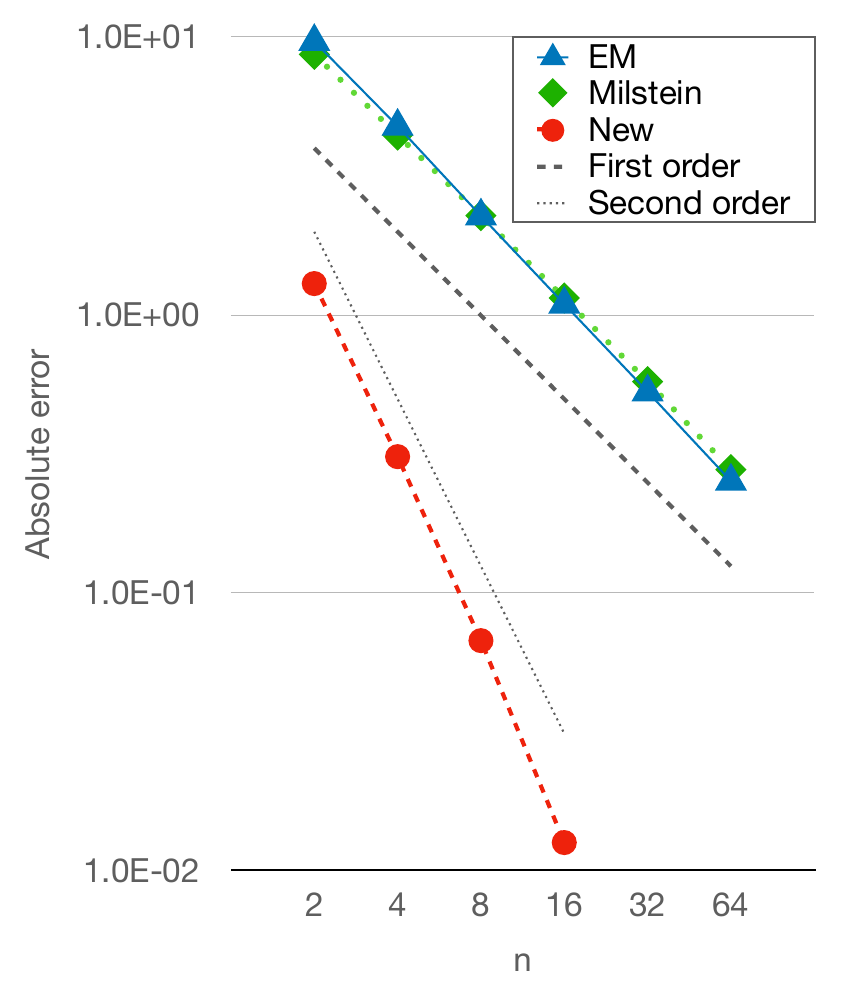}
   \caption{Weak convergence ($\sigma = 0.8$).} 
    \label{fig:asian_bs_conv_2}
\end{subfigure}  
     \caption{Asian call option pricing (BS model). }
     \label{fig:asian_bs}
\end{figure}

We consider the following 2-dimensional hypo-elliptic SDE:
\begin{align} \label{eq:asian_bs}
\begin{aligned} 
d S_t & = r S_t dt + \sigma S_t  d B_t^1; \\ 
d A_t & = S_t dt, 
\end{aligned}
\end{align}
with the parameters $r > 0$ and $\sigma > 0$. $S_t$ and $A_t$ represent the price of the underlying asset and the time accumulation of the asset value, respectively. Then, the price of \emph{Asian call option} with the strike price $K$ and the maturity $T$ is formulated as 
$
  D \times \E \bigl[ \varphi_K (A_T / T) \bigr], 
$
where the text function is $\varphi_K (x) = \max (x - K, 0)$ for $x \in \mathbb{R}^N$, and $D = e^{-rT}$ is the discount factor. We note that for the SDE (\ref{eq:asian_bs}) the commutative condition holds and thus the Milstein scheme coincides with the truncated Milstein scheme. We consider the following setting for parameter values: $r = 0.1$, $\sigma \in 
\{0.4, 0.8 \}$, $T = 1.0$ and $(S_0, A_0) = (100.0, 0.0)$. We compute the benchmark value by the standard Monte-Carlo (MC) method using the EM scheme with the number of paths $M =10^7$ and discretization $n = 2^{10}$. We compare the performances of three numerical schemes, EM, Milstein and the extended Milstein scheme via the following MC estimates: for a number of MC paths $M$ and discretisation $n$, 
\begin{align*}
f^{w}_K  (M, n)  = \tfrac{1}{M} \sum_{j = 1}^M D \times 
\varphi_K \bigl( \widetilde{A}^{w, n, [j]}_T / T \bigr), \qquad w \in \{\mrm{EM}, 
\mrm{Mil}, \mrm{New} \},  
\end{align*} 
where $\widetilde{A}^{\mrm{EM}, n, [j]}_T, \, \widetilde{A}^{\mrm{Mil}, n, [j]}_T, \, \widetilde{A}^{\mrm{New},n, [j]}_T$ are the $j$-th trajectory of the EM scheme, the truncated Milstein scheme and the extended Milstein scheme (\ref{eq:ext-Mil}) applied to the model (\ref{eq:asian_bs}), respectively.  
In Figures \ref{fig:asian_bs_err_1} and \ref{fig:asian_bs_err_2}, 
we plot 
$$
\mathscr{E}_K^w (M, n ) \equiv (\mrm{Benchmark \, value})_K 
- f^w_K (M, n), \qquad  w \in \{ \mrm{EM}, \mrm{Mil}, \mrm{New} \}, 
$$ 
with $M = 10^6$ for $n = 4, 8, 16$ and $K = 10, 20, \ldots, 200$, where we applied Quasi-Monte-Carlo (QMC) for computing $f^w_K (M, n)$. We observe that the proposed scheme provides very accurate estimates for all the strike prices and outperforms the EM/Milstein schemes. Also, Figures \ref{fig:asian_bs_conv_1} and \ref{fig:asian_bs_conv_2} show 
\begin{align} 
\sup_{K \in \{10, 20, \ldots, 200\}} \bigl| \mathscr{E}_K^w (10^6, n) \bigr|,  \qquad  w \in \{ \mrm{EM}, \mrm{Mil}, \mrm{New} \},  
\end{align}
for various numbers of discretisation $n$. It is noteworthy that there is no significant difference between the performance of the EM and Milstein schemes. 
\subsection{Asian digital option pricing with Heston model}
We consider the following 3-dimensional SDE: 
\begin{align}
\begin{aligned} \label{eq:asian_heston}
d S_t & = \sqrt{\sigma_t} S_t dB_t^1; \\[0.1cm]  
d \sigma_t & = \alpha (\theta - \sigma_t) dt 
+  \nu \sqrt{\sigma_t} ( \rho dB_t^1 + \sqrt{1 - \rho^2} dB_t^2);  \\[0.1cm] 
d A_t & = S_t dt,  
\end{aligned}
\end{align} 
with parameters $\alpha, \nu, \theta > 0$ and $\rho \in [-1, 1]$ satisfying $2 \alpha \theta > \nu^2$ so that the process $\{\sigma_t \}_{t \geq 0}$ is strictly positive. The pair $\{S_t, \sigma_t\}_{t \ge 0}$ represents the process of the underlying asset and its volatility, and is a famous stochastic volatility model called \emph{Heston model} in the context of financial mathematics. We then estimate the price of of Asian digital option with the coupon $\mrm{Cpn} > 0$, the strike price $K > 0$ and the maturity $T > 0$, given as: 
\begin{align} \label{eq:digtal_option}
\mrm{Cpn} \times \E \bigl[ H_K (A_T / T) \bigr], 
\end{align} 
where $H_K (x) \equiv \mathbf{1}_{x \geq K}$. Note that the test function is non-smooth. Since (\ref{eq:digtal_option}) does not admit a closed-form solution, we apply the quasi-Monte Carlo method accompanied by some discretisation scheme to approximate the quantity. In particular, the commutative condition does not hold for the model (\ref{eq:asian_heston}), and then the standard Milstein scheme is intractable due to the presence of L\'evy area. We thus compare the performance of these three numerical schemes via the following Monte-Carlo estimates:
$$
f^{w}_K  (M, n)  = \tfrac{1}{M} \sum_{j = 1}^M \mrm{Cpn} \times 
H_K \bigl( \widetilde{A}^{w, n, [j]}_T / T \bigr), \qquad w \in \{\mrm{EM}, 
\mrm{TMil}, \mrm{New} \}  
$$ 
with the number of trajectories $M$ and the number of discretisation $n$, where $\widetilde{A}^{\mrm{EM}, n, [j]}_T, \, \widetilde{A}^{\mrm{TMil}, n, [j]}_T, \, \widetilde{A}^{\mrm{New},n, [j]}_T$ are the $j$-th trajectory of the EM scheme, the truncated Milstein scheme and the extended Milstein scheme (\ref{eq:ext-Mil}) applied to the model (\ref{eq:asian_heston}), respectively. We set the parameter values as: 
$\alpha = 2.0$, $\theta = 0.09$, $\nu = 0.1$, $\rho = 0.7$, $T = 1.0$, $
(S_0, \sigma_0, A_0)= (100, 0.09, 0.0)$ and $\mrm{Cpn} = 100$. We estimate the benchmark value by applying the standard Monte-Carlo method with the EM scheme as:
$(\mrm{Benchmark value})_K$ = $f^{\mrm{EM}}_K (10^7, 2^{11})$. 
\begin{figure}[h]
     \centering
     \begin{subfigure}[b]{0.55\textwidth}
         \centering
         \includegraphics[width=\textwidth]{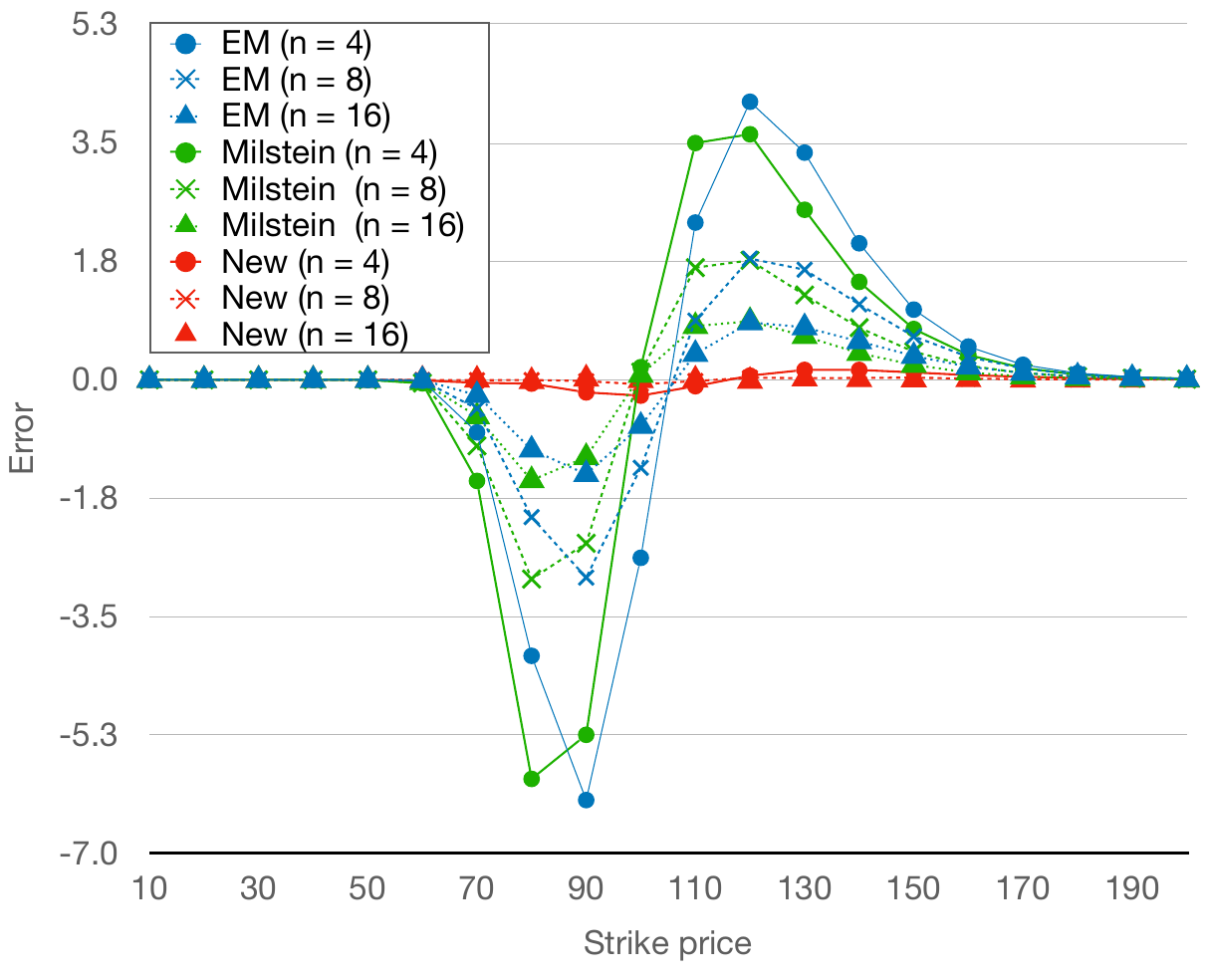}
        \caption{Error of estimation.} %
         \label{fig:heston_err}
     \end{subfigure}
     \hspace{0.5cm}
     \begin{subfigure}[b]{0.38\textwidth}
         \centering
         \includegraphics[width=\textwidth]{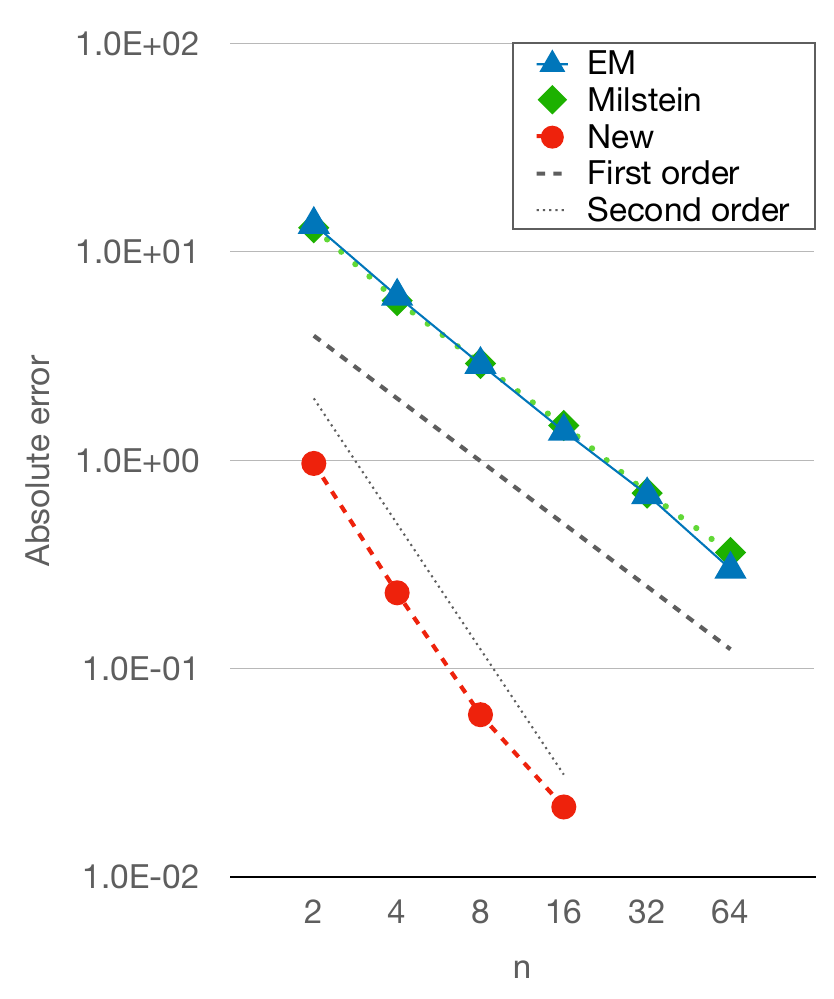}
        \caption{Weak convergence.} 
         \label{fig:heston_conv}
     \end{subfigure} 
        \caption{Asian digital option pricing (Heston model) with three numerical schemes. }
        \label{fig:mle}
\end{figure}
In Figure \ref{fig:heston_err}, we plot 
$$
\mathscr{E}_K^w (M, n ) \equiv (\mrm{Benchmark \, value})_K 
- f^w_K (M, n), \qquad  w \in \{ \mrm{EM}, \mrm{TMil}, \mrm{New} \}, 
$$ 
with $M = 10^6$ for $n = 4, 8, 16$ and $K = 10, 20, \ldots, 200$, where we applied Quasi-Monte-Carlo (QMC) for computing $f^w_K (M, n)$. Also, Figure \ref{fig:heston_conv} indicates 
\begin{align} 
\sup_{K \in \{10, 20, \ldots, 200\}} \bigl| \mathscr{E}_K^w (10^6, n) \bigr|,  \qquad  w \in \{ \mrm{EM}, \mrm{TMil}, \mrm{New} \},  
\end{align}
for various numbers of discretization $n$. As we observed in the previous numerical experiments on the Asian BS model, the proposed scheme again provides quite accurate estimates compared with the EM scheme and the truncated Milstein scheme. Notably, from Figure \ref{fig:heston_conv}, the bias induced by the proposed scheme is less than a tenth of that by the EM scheme or the truncated Milstein scheme, for $n \in \{2, 4, 8, 16\}$ and achieves nearly the second order weak convergence, which is expected from the argument in Section \ref{sec:small_diff} for small diffusions.  
\section{Proof of Theorem \ref{thm:wa_mil}} \label{sec:pf_main}
Most of our proof relies on the argument in \citep{bally96} that shows Theorem \ref{thm:wa_em}, i.e., weak convergence of EM scheme under non-smooth test functions. The strategy is naturally applied to the case of extended Milstein scheme (\ref{eq:ext-Mil}), and our main focus is to derive an explicit form of the leading error term of weak approximation by the extended Milstein scheme, which is associated with the local error expansion by the scheme (Proposition \ref{prop:loc_mil}). The technical results introduced later (Lemmas \ref{lemma:ibp_Milstein}, \ref{lemma:ibp_Milstein_2} and \ref{lemma:est_last_step}) are also shown by following the proofs in \citep{bally96} and adjusting the argument for EM scheme with our proposed scheme (\ref{eq:ext-Mil}). Thus, the detailed proof of these results will be omitted.  

We introduce some notation for the proof. Let $T >0$, $h = T/n, \, n \ge 1$ and we assume $f \in \mathscr{B}_b (\mathbb{R}^N)$. We introduce families of Markov operators $\{P_t\}_{t \ge 0}$ and $\{Q_t \}_{t \ge 0}$ as:
\begin{align}
P_t f (x) = \E [f (X_t^x)], \quad 
Q_t f (x) = \E [f ( \bar{X}_t^x)],  \quad t > 0, \, x \in \mathbb{R}^N, 
\end{align} 
with $\bar{X}_t^x$ being the one step extended Milstein scheme specified as: 
\begin{align} \label{eq:one_step_new}
\bar{X}_t^x = x + \sum_{j = 0}^d \sigma_j (x) B_t^j 
+ \sum_{0 \le j_1, j_2 \le d} L_{j_1} \sigma_{j_2} (x) \times 
\tfrac{1}{2} (B_t^{j_1} B_t^{j_2} - t \mathbf{1}_{j_1 = j_2 \neq 0 }),
\end{align}
where we recall $\sigma_0 = b$ and $L_0 = \mathcal{L}$ given in (\ref{eq:L}). We note that $k$-times application of the Markov operator $Q_{h}$ gives: 
\begin{align} 
 (Q_h)^k \varphi (x) = \E [\varphi (\new{n}{x}{kh})], \quad \varphi \in \mathscr{B}_b (\mathbb{R}^N), 
\end{align}
where $\new{n}{x}{kh}$ is the Markov chain of the one step extended Milstein scheme defined as (\ref{eq:ext-Mil}). First, we have the following result associated with the local weak approximation by the extended Milstein scheme (\ref{eq:ext-Mil}):
\begin{prop} \label{prop:loc_mil}
Let $x \in \mathbb{R}^N$ and $h = T/n, \, T > 0, \, n \in \mathbb{N}$. We have  
\begin{align} \label{eq:loc_mil} 
(Q_h)^k (P_h - Q_h) 
P_{T - (k+1) h} f (x)
=  h^2 \, 
\E \Bigl[ \Phi_{3} \bigl( kh, \new{n}{x}{kh}
\bigr)
\Bigr]  + \widetilde{\Phi}^{(n)}_{k} (h, x), \quad  0 \le k \le n-2,  
\end{align} 
where $\Phi_3$ is defined as (\ref{eq:Phi_1}), and $\widetilde{\Phi}^{(n)}_{k} (h, x)$ is determined from the sum of terms given in the form of
\begin{align} 
\label{eq:err}
& \mathbb{E} \left[ e_\alpha (\new{n}{x}{kh}) 
\int_{kh}^{(k+1)h} \int_{kh}^{s_1} \int_{kh}^{s_2}
g_\alpha ( \new{n}{x}{s_3} ) \partial_\alpha u (s_3, \new{n}{x}{s_3}) 
ds_3 ds_2 ds_1  \right], 
\end{align}
where $\alpha \in \{1, \ldots, N\}^\ell$ with $\ell \le 6$, and the functions $e_\alpha, \, g_\alpha$ are given as products of partial derivatives of coefficients $\sigma_j^i, \, 0 \le j \le d, \, 1 \le i \le N$.   
\end{prop}
The proof of Proposition \ref{prop:loc_mil} is provided in Appendix \ref{appendix:loc_mil}. Proposition \ref{prop:loc_mil} gives the global weak error expansion for the proposed scheme (\ref{eq:ext-Mil}) in the following form: 
\begin{align}
& P_T f (x)  - (Q_h)^n f (x) 
= \sum_{k = 0}^{n-1} (Q_h)^k (P_h - Q_h) P_{T - (k+1)h}f(x) \nonumber \\
& \qquad  = h^2 \sum_{k = 0}^{n-2}  
\E 
\bigl[
  \Phi_{3} ( kh, \new{n}{x}{kh} )  
\bigr]   
+ \sum_{k = 0}^{n-2} \widetilde{\Phi}_{k}^{(n)} (h, x) 
+ (Q_{h})^{n-1} (P_{h} - Q_{h}) f(x) 
\label{eq:global_err_0}  \\
& \qquad  = h \int_0^T \E \bigl[ \Phi_{3} ( s, X_s^x ) \bigr] ds 
+ \sum_{i = 1}^4 \mathcal{E}^f_{i,n} (h, x),  
\label{eq:global_err}
\end{align}
where we have set: 
\begin{align*}
\mathcal{E}^f_{1,n} (h, x) 
& \equiv 
h^2 \sum_{k = 0}^{n-2}  
\E \bigl[ \Phi_{3} (kh, X_{kh}^x) \bigr] 
- h \int_0^T \mathbb{E} \bigl[ \Phi_{3} ( s, X_s^x ) \bigr] ds; \\ 
\mathcal{E}^f_{2,n} (h, x) 
& \equiv h^2 \sum_{k = 0}^{n-2}  
\E \bigl[
  \Phi_{3} ( kh, \new{n}{x}{kh} )  
\bigr] 
- h^2 \sum_{k = 0}^{n-2}  
\E 
\bigl[ \Phi_{3} ( kh, X_{kh}^x ) \bigr];   \\ 
\mathcal{E}^f_{3,n} (h, x) 
& \equiv \sum_{k = 0}^{n-2} \widetilde{\Phi}_{k}^{(n)} (h, x), 
\quad 
\mathcal{E}^f_{4,n} (h, x) \equiv 
(Q_{h})^{n-1} (P_{h} - Q_{h}) f(x). 
\end{align*}
To show the upper bounds of the terms $\mathcal{E}_{i,n}^f (h, x), \; 1 \le k \le 4$, we introduce the following three lemmas (corresponding results based upon the EM scheme are found in \cite{bally96} as Lemmas 4.1, 4.2 and 4.3, respectively): 
%
%
\begin{lemma} \label{lemma:ibp_Milstein}
For any multi-index $\alpha \in \{1, \ldots, N\}^\ell, \, \ell \in \mathbb{N}$ and any smooth function $g : \mathbb{R}^N \to \mathbb{R}$ with the polynomial growth, there exist a non-decreasing function $K(\cdot)$ and constants $q, Q > 0$ independent of $n, T$ such that 
\begin{align} \label{eq:bd_deriv_diff}
\bigl| \E \bigl[ g (X_s^x ) \partial_\alpha P_{T-s} f (X_s^x) \bigr] \bigr| \leq K(T) \frac{\| f \|_\infty}{T^q} ( 1 +  |x|^Q)
\end{align} 
for all $s \in [0, T]$ and
\begin{align} \label{eq:bd_deriv_mil}
\bigl| \mathbb{E} \bigl[ g (\new{n}{x}{s} ) \partial_\alpha 
P_{T-s} f (\new{n}{x}{s}) \bigr] \bigr| \leq K(T) \frac{\| f \|_\infty}{T^q} ( 1 +  |x|^Q)
\end{align}
for all $s \in [0, T- T/n]$. 
\end{lemma}


\begin{lemma} \label{lemma:ibp_Milstein_2}
Let $\alpha \in \{1, \ldots, N \}^{l_1}, \, \beta \in  \{1, \ldots, N \}^{l_2}$, $l_1, l_2 \in \mathbb{N}$ and 
$g, \widetilde{g} : \mathbb{R}^N \to \mathbb{R}$ be smooth functions with polynomial growth. Set $\varphi_s : \mathbb{R}^N \to \mathbb{R}, \, s \in [0, T]$ as:
$ y \mapsto \varphi_s (y) \equiv  \widetilde{g} (y) \partial_\beta P_{T- s} f (y). $
Then, there exist a non-decreasing function $K(\cdot)$ and constants $q, Q > 0$ independent of $n, T$ such that for all $s \in [0, T-T/n]$ and for all $t \in [0, s - T/n]$, 
\begin{align}
\bigl| 
\mathbb{E} [ g (\bar{X}_t^{n , x}) \partial_\alpha P_{s - t} \varphi_s (\bar{X}_t^{n , x}) ] 
\bigr| 
\le K(T) \frac{\| f \|_\infty}{T^q} (1 + |x|^Q). 
\end{align}
\end{lemma}

\begin{lemma} \label{lemma:est_last_step}
There exist a non-decreasing function $K(\cdot)$ and a constant $Q > 0$ independent of $n, T$ such that 
\begin{align} 
 \bigl| \mathcal{E}_{4,n}^f (h,x)  \bigr| 
 \leq K(T) \frac{\| f \|_\infty }{n^2} (1 + |x|^Q). 
\end{align}
\end{lemma}

\noindent 
(\textit{Proof of Lemmas \ref{lemma:ibp_Milstein}, \ref{lemma:ibp_Milstein_2} and \ref{lemma:est_last_step}.})
These are shown via the similar argument used to prove Lemmas 4.1, 4.2, 4.3 in \citep{bally96} by replacing the Euler-Maruyama scheme with the extended Milstein scheme (\ref{eq:ext-Mil}), and this replacement does not require special treatments. We thus omit the detailed proof.  
$\Box$
\\ 

We now return to the proof of Theorem \ref{thm:wa_mil}. Since the upper bound for the term  $\mathcal{E}_{4,n}^f (h, x)$ is provided in Lemma \ref{lemma:est_last_step}, we will study the other three terms. For the first term, we have $\mathcal{E}^f_{1,n} (h, x)  
= \mathcal{E}^{f, (\mrm{I})}_{1,n} (h, x)  
+ \mathcal{E}^{f,(\mrm{II})}_{1,n} (h, x)$, 
with  
\begin{align*}
    \mathcal{E}^{f, (\mrm{I})}_{1,n} (h, x)  
    & \equiv 
    - h^2 \mathbb{E} 
    \bigl[
        \Phi_{3} ( (n-1)h, X_{(n-1)h}^x)  
    \bigr], \\
     \mathcal{E}^{f, (\mrm{II})}_{1,n} (h, x)  
    & \equiv  
    h^2 \sum_{k = 0}^{n-1}  
  \mathbb{E} \bigl[
      \Phi_{3} ( k h, X_{kh}^x )  
    \bigr] 
   - h \int_0^T \mathbb{E} \bigl[ \Phi_{3} ( s, X_s^x ) \bigr] ds.  
\end{align*}
From the definition of $\Phi_3$ and the bound (\ref{eq:bd_deriv_diff}) in Lemma \ref{lemma:ibp_Milstein}, we have: 
\begin{align}
\bigl| \mathcal{E}^{f, (\mrm{I})}_{1,n} (h, x) \bigr| \le 
K(T) \frac{\| f \|_\infty}{n^2} (1 + |x|^Q). 
\end{align} 
We also have: 
\begin{align*}
\Bigl| \frac{T}{n} \sum_{k = 0}^{n-1}  
  \mathbb{E} \bigl[
      \Phi_{1} ( kT/n, X_{kT/n}^x )  
    \bigr] 
   - \int_0^T \mathbb{E} \bigl[ \Phi_{1} ( s, X_s^x ) \bigr] ds   \Bigr| \leq K(T) \frac{\| f \|_\infty}{n} (1 + |x|^Q),
\end{align*}
where we applied It\^o formula as
\begin{align*} 
\int_{k h}^{(k+1)h} \E 
\bigl[ \Phi_3 (s, X_s^x)  \bigr] ds 
& = h \E 
\bigl[ \Phi_3  (k h, X_{kh}^x) \bigr]
+ \int_{k h}^{(k+1)h} \int_{kh}^s 
\E \bigl[ \partial_t \Phi_3 (s_1, X_{s_1}^x)  \bigr] ds_1 ds   \\[0.2cm]
& \quad 
+ \int_{kh}^{(k+1)h} \int_{k h}^s 
\mathbb{E} \bigl[ \mathcal{L} \Phi_3 (s_1, X_{s_1}^x)  \bigr] ds_1 ds,  
\quad k = 0, \ldots, n-1,  
\end{align*}
and the last two terms in the right-hand side are bounded by $K(T) \tfrac{\| f \|_\infty}{n^2} (1 + |x|^Q)$ due to (\ref{eq:bd_deriv_diff}). We thus obtain the upper bound 
$\bigl| \mathcal{E}_{1, n}^f (h, x) \bigr| 
\le K(T) \tfrac{\|f \|_\infty}{n^2} (1 + |x|^Q)$. 

Next, we study the second error term $\mathcal{E}_{2, n}^f (h, x)$. The term writes $\mathcal{E}_{2, n}^f (h, x) 
= \mathcal{E}_{2, n}^{f, (\mathrm{I})} (h, x) + \mathcal{E}_{2, n}^{f, (\mathrm{II})} (h, x)$ with 
\begin{align*} 
\mathcal{E}_{2, n}^{f, (\mathrm{I})} (h, x) 
& \equiv  h^2 \sum_{0 \le k \le [n/2]} 
\Bigl\{ \E \bigl[
\Phi_{3} ( kh, \new{n}{x}{kh} )  
\bigr] - 
\E \bigl[
\Phi_{3} ( kh, X_{kh}^x )  \bigr] 
\Bigr\}; \\   
\mathcal{E}_{2, n}^{f, (\mathrm{II})} (h, x) 
& \equiv  h^2 \sum_{[n/2] + 1 \le k \le n-2} 
\Bigl\{ 
\E \bigl[
\Phi_{3} ( kh, \new{n}{x}{kh} )  
\bigr] - 
\E \bigl[
\Phi_{3} ( kh, X_{kh}^x )  \bigr] 
\Bigr\}.  
\end{align*}
We consider the term $\mathcal{E}_{2, n}^{f, (\mathrm{I})} (h, x)$. First, applying (\ref{eq:bd_deriv_mil})  and Lemma \ref{lemma:est_last_step} to (\ref{eq:global_err_0}), we have that:  for any $\varphi \in \mathscr{B}_b(\mathbb{R}^N)$, there exist a non-decreasing function $K (\cdot)$ and constants $q, Q > 0$ such that for all $s \in [0, T]$, 
\begin{align} \label{eq:first_order_bd}
\bigl|  \mathbb{E} [\varphi (X_s^x)] 
- \mathbb{E} [\varphi (\bar{X}_s^{n, x})]
\bigr| \le \tfrac{K(T)}{T^q}
\tfrac{\| \varphi \|_\infty}{n} (1 + |x|^Q).  
\end{align}
Using the bound (\ref{eq:first_order_bd}) with $\varphi = \Phi_3 (kh, \cdot)$, we get
\begin{align} \label{eq:err_2_1}
\bigl| \mathcal{E}_{2, n}^{f, (\mathrm{I})} (h, x)  \bigr|
& \le h^2 \sum_{0 \le k \le [\tfrac{n}{2}]} 
\tfrac{K_1 (T)}{T^{q_1}} \cdot \tfrac{\| \Phi_3 (kh, \cdot) \|_\infty}{n} (1 + |x|^{Q_1}) 
\le  \tfrac{K_2 (T)}{T^{q_2}} \cdot \tfrac{\| f  \|_\infty}{n^2} (1 + |x|^{Q_2}),  
\end{align}
for some non-decreasing functions $K_1, K_2$ and constants $q_1, q_2, Q_1, Q_2 > 0$ independent of $n, T$. We have also used the estimate (4.9) in \citep{bally96} to bound $\| \Phi_3 (kh, \cdot) \|_\infty$ in the last inequality. We next study the term $\mathcal{E}_{2, n}^{f, (\mathrm{II})} (h, x)$. We apply the error expression (\ref{eq:global_err_0}) to the term $\E \bigl[
\Phi_{3} ( kh, \new{n}{x}{kh} )  
\bigr] - \E \bigl[\Phi_{3} ( kh, X_{kh}^x ) \bigr]$
by replacing the test function $f (\cdot)$ with $\Phi_3 (kh, \cdot)$, and then the term is ultimately bounded by $\tfrac{K(T)}{T^q} \tfrac{\| f \|_\infty}{n}(1 + |x|^Q)$ due to Lemma \ref{lemma:ibp_Milstein_2} with some non-decreasing function $K(\cdot)$ and constants $q, Q > 0$ independent of $n , T$. This immediately leads to that the term $\bigl| \mathcal{E}_{2, n}^{f, (\mathrm{II})} (h, x) \bigr|$ is bounded as the right-hand side of (\ref{eq:err_2_1}). 
Thus, we have
\begin{align*} 
\bigl| \mathcal{E}_{2, n}^f (h , x) \bigr| 
& \leq 
h^2 \sum_{k = 0}^{n-2} 
\Bigl| \mathbb{E} 
\bigl[
\Phi_{3} ( kh, \new{n}{x}{kh} )  
\bigr] 
- \mathbb{E} 
\bigl[
\Phi_{3} ( kh, X_{kh}^x )  
\bigr]    \Bigr| 
\leq  K(T) \frac{\| f \|_\infty}{n^2} ( 1 + |x|^Q). 
\end{align*}
%

Finally, we study the term $\mathcal{E}_{3, n}^f (h, x)$. Notice that each term in  $\widetilde{\Phi}_k^{(n)} (x)$ contains triple time integrals inside the expectation. Thus, Lemma \ref{lemma:ibp_Milstein} yields
\begin{align}
 \Bigl| \mathcal{E}_{3, n}^f (h, x) \Bigr| 
 \leq \sum_{k = 0}^{n-1} \bigl| \widetilde{\Phi}_k^{(n)} (x) \bigr| 
 \leq K(T) \frac{\| f \|_\infty}{n^2} (1 + |x|^Q),
\end{align}
and we conclude. 
\section{Conclusions}
We have proposed a straightforward and effective first order weak scheme (\ref{eq:ext-Mil}) for diffusion processes and compared its efficacy with other popular first order schemes, specifically, Euler-Maruyama/Misltein/truncated Milstein schemes. The proposed scheme is always explicit and can be simulated with Gaussian random variables, with the same number of random variables required by the Euler-Maruyama scheme or the truncated Milstein scheme.  We have shown that the proposed scheme achieves the first order weak convergence, but its leading order error term involves fewer terms than other first order schemes. In particular, when the diffusion coefficients contain a small parameter $\varepsilon$, the discretization bias can be significantly reduced with the aid of $\varepsilon$, though the effect does not necessarily appear for EM and (truncated) Milstein schemes. We then carried out numerical experiments of Asian option pricing, which shows that the new scheme provides much superior accuracy to the other first order schemes and behaves nearly as the second order scheme under the setting of small diffusions. Notably, we have also observed that there are no significant differences in the performance between the EM and Milstein schemes in the experiment. 

The application of the proposed time-discretization scheme is not limited to Monte-Carlo estimate for the expectation of the law of diffusions. Indeed, developing simple and effective discretization is useful in a much wider context that requires an approximate sampling of diffusion processes, e.g., parameter estimation of SDEs or filtering of diffusion processes when only partial coordinates are observed. For instance, \cite{igu22} recently emphasized that the use of time-discretization with accurate weak approximation leads to efficient parameter estimation of diffusions when the so-called \emph{Data augmentation} approach \citep{DA13} is required to conduct Bayesian inference from low frequency observations. Due to the simple definition and its effective weak approximation, our proposed scheme can be incorporated into many computational/statistical methodologies (e.g. \citep{manifold22, AMLMC24}). 

\section*{Acknowledgements}
Yuga Iguchi is supported by the Additional Funding Programme for Mathematical Sciences, delivered by EPSRC (EP/V521917/1) and the Heilbronn Institute for Mathematical Research. Toshihiro Yamada was supported by JST PRESTO (JPMJPR2029), Japan. We thank Prof. Shigeo Kusuoka (University of Tokyo) for grateful advice for the method of the paper.
\appendix
%
\section{Proof of Lemma \ref{lemm:Phi_em}} \label{appendix:Phi_em} 
Making use of the generator $\mathcal{L}$ given in (\ref{eq:L}) and $\partial_t u = - \mathcal{L} u $, we have from (\ref{eq:err_em_original}) that
\begin{align} \label{eq:em_err}
\begin{aligned}
\Phi^{\mathrm{EM}} (t,x) 
& =  - \frac{1}{2} \sum_{i, j = 1}^N b^i (x) b^j (x) \partial_{i j} u (t, x) - \frac{1}{2} \sum_{i, j, k  = 1}^N b^i (x)  a^{jk}(x) \partial_{ijk} u (t,x) \\ 
& \quad - \frac{1}{8} \sum_{i,j,k,l = 1}^N a^{ij}(x)  a^{kl}(x) \partial_{ijkl} u (t, x) 
+ \frac{1}{2} \mathcal{L}^2 u (t, x). 
\end{aligned}
\end{align} 
Since 
   \begin{align*}
      & \mathcal{L}^2 u (t,x) 
       = \sum_{i, j = 1}^N b^i (x) 
      \bigl\{ \partial_i b^j (x) \partial_j u (t, x)  +  b^j(x) \partial_{ij} u (t, x)  \bigr\} \\  
      & \quad 
      + \frac{1}{2} \sum_{i, j, k = 1}^N b^i (x) \bigl\{ \partial_i a^{jk}(x) \partial_{jk} u (t,x)  + a^{jk}(x) \partial_{ijk} u (t,x) \bigr\} \\
      & \quad + \frac{1}{2} \sum_{i,j,k = 1}^N 
     a^{ij}(x)  
      \bigl\{  \partial_{ij} b^k(x) \partial_{k} u(t,x) 
      + \partial_i b^k (x) \partial_{jk} u (t,x) 
        + \partial_j b^k (x) \partial_{ik} u (t,x) 
        + b^k(x) \partial_{ijk} u (t,x)   \bigr\}  \\
      & \quad + \frac{1}{4} \sum_{i,j,k,l  = 1}^N 
        a^{ij}(x)  
         \bigl\{  \partial_{ij} a^{kl} (x) \partial_{kl} u(t,x) 
         + \partial_i a^{kl} (x) \partial_{jkl} u (t,x) 
           + \partial_j a^{kl} (x) \partial_{ikl} u (t,x) 
           + a^{kl} (x) \partial_{ijkl} u (t,x)   \bigr\},  
\end{align*}
it follows that
\begin{align} \label{eq:Phi_em_2}
\Phi^{\mathrm{EM}} (t,x)
& =
\frac{1}{2} \sum_{i, j = 1}^N b^i (x) 
\partial_i b^j (x) \partial_j u (t, x)  
+  \frac{1}{4} \sum_{i, j, k = 1}^N b^i (x)  
\partial_i a^{jk}(x) \partial_{jk} u (t,x) \nonumber  \\
& \quad 
+  \frac{1}{4} \sum_{i,j,k = 1}^N 
a^{ij}(x)  \partial_{ij} b^k(x) \partial_{k} u(t,x) 
+ \frac{1}{2} \sum_{i,j,k = 1}^N 
a^{ij}(x) \partial_i b^k (x) \partial_{jk} u (t,x)  
\nonumber  \\
& \quad +  \frac{1}{8} \sum_{i,j,k,l  = 1}^N   a^{ij}(x)  
 \partial_{ij} a^{kl} (x) \partial_{kl} u(t,x) 
+ \frac{1}{4} \sum_{i,j,k,l  = 1}^N 
a^{ij}(x) \partial_i a^{kl} (x) \partial_{jkl} u (t,x). 
\end{align} 
Substituting 
\begin{align*}
    & \partial_{i} a^{jk}(x) 
    = \sum_{m = 1}^d \Bigl\{  \partial_{i} \sigma_m^j (x) \, \sigma_m^k (x) 
    + \sigma_m^j (x) \partial_{i} \sigma_m^k (x)    \Bigr\}, \\ 
    & \partial_{ij} a^{kl}(x) 
    = \sum_{m = 1}^d \Bigl\{ 
      \partial_i \sigma_m^k (x) \, \partial_j \sigma_m^l (x) 
    + \partial_j \sigma_m^k (x) \, \partial_i \sigma_m^l (x) 
    + \partial_{ij} \sigma_m^k (x) \,   \sigma_m^l (x) 
    + \sigma_m^k (x) \,  \partial_{ij} \sigma_m^l (x) 
   \Bigr\} 
\end{align*}
into the right hand side of the equation (\ref{eq:Phi_em_2}), we have: 
\begin{align*}
\Phi^{\mathrm{EM}} (t,x)
& = \tfrac{1}{2} \sum_{i = 1}^N L_0 b^i (x) \partial_i u (t,x)
+ \tfrac{1}{2} \sum_{i,j = 1}^N \sum_{m = 1}^d 
\sigma_{m}^i (x) L_{m} b^j (x)\partial_{ij} u (t,x) 
+ \tfrac{1}{2} \sum_{i,j = 1}^N \sum_{m = 1}^d  \sigma_m^i(x) L_0 \sigma_m^j (x) \partial_{ij} u (t,x) \nonumber \\ 
& \quad  + \tfrac{1}{4} \sum_{i,j=1}^N \sum_{m_1, m_2}^d 
L_{m_1} \sigma_{m_2}^i (x) L_{m_1} \sigma_{m_2}^j (x) 
\partial_{ij} u(t,x)   
+  \tfrac{1}{2} \sum_{i, j, k = 1}^N \sum_{m_1, m_2}^d L_{m_1} \sigma_{m_2}^{i} (x) \sigma_{m_1}^j (x) 
\sigma_{m_2}^k (x) \partial_{ijk} u(t,x) \\
& = \Phi_1 (t,x) + \Phi_2 (t,x) + \Phi_3 (t,x), 
\end{align*}
where the last line is deduced from:
\begin{align*}
\tfrac{1}{4} \sum_{i,j=1}^N \sum_{m_1, m_2}^d 
L_{m_1} \sigma_{m_2}^i (x) L_{m_1} \sigma_{m_2}^j (x) 
\partial_{ij} u(t,x) 
& = \tfrac{1}{8} \sum_{i,j=1}^N \sum_{m_1, m_2}^d 
L_{m_1} \sigma_{m_2}^i (x) \Bigl\{
 L_{m_1} \sigma_{m_2}^j (x) 
+  L_{m_2} \sigma_{m_1}^j (x)  
\Bigr\}
\partial_{ij} u(t,x)  \\
& \quad 
+ \tfrac{1}{8} \sum_{i,j=1}^N \sum_{m_1, m_2}^d 
L_{m_1} \sigma_{m_2}^i (x) \Bigl\{
L_{m_1} \sigma_{m_2}^j (x) 
-  L_{m_2} \sigma_{m_1}^j (x)  
\Bigr\} \partial_{ij} u(t,x),  
\end{align*} 
and we conclude. 
\section{Proof of Proposition \ref{prop:loc_mil}} \label{appendix:loc_mil}
It holds that for $k = 0,1, \ldots, n-2$, 
\begin{align} \label{eq:local_err}
(Q_{h})^k (P_{h} - Q_{h}) P_{T- (k+1)h} f(x)
= \E \bigl[ P_{T-k h} f ( \new{n}{x}{kh} ) \bigr] 
- \mathbb{E} \Bigl[ P_{T-(k + 1)h} 
f \bigl( \new{n}{x}{(k + 1)h} \bigr) \Bigr]  
= - \E \Bigl[ F_k^{(n)} \bigl( \new{n}{x}{kh} \bigr) \Bigr], 
\end{align}
where we have set: 
\begin{align*}
F_k^{(n)} ( z ) 
\equiv \mathbb{E} \Bigl[ 
u ( (k+1) h,  \bar{X}_h^z ) - u ( kh,  z) \Bigr], \ \ z \in \mathbb{R}^N,
\end{align*}
with $\bar{X}_h^z$ being the one step extended Milstein scheme (\ref{eq:one_step_new}) given the previous state $z$. We now consider the expansion of $F_k^{(n), x} ( z )$ via iterative application of It\^o's formula to $u ( (k+1)h,  \mil{z}{h})$. Application of It\^o's formula yields:
\begin{align*}
F_k^{(n)} ( z )   
= T_1  + T_2 + T_3,
\end{align*}
where we have defined: 
\begin{gather*}
T_1 = \E \Bigl[ \int_0^h \partial_t u (kh + s, \bar{X}^z_s ) ds \Bigr], \qquad 
T_2 = \sum_{i = 1}^N \E \Bigl[ \int_0^h \partial_i u (kh + s, \bar{X}^z_s ) d \bar{X}^{z, i}_s \Bigr], \\ 
T_3 = \tfrac{1}{2} 
\sum_{i, j  = 1}^N \E \Bigl[ \int_0^h \partial_{ij} u (kh + s, \bar{X}^z_s ) d \langle \bar{X}^{z, i},   \bar{X}^{z, j} \rangle_s  \Bigr]. 
\end{gather*} 
We have that for $t >0$ and $x \in \mathbb{R}^N$, 
\begin{align*}
& d \bar{X}^{x, i}_{t}
= \sum_{k=0}^d \sigma_k^{i} (x) d B_t^k 
+ \tfrac{1}{2} \sum_{0 \le k_1, k_2 \le d} 
 \bigl\{ L_{k_1} \sigma_{k_2}^{i} (x) + L_{k_2} \sigma_{k_1}^{i} (x) \bigr\} B_t^{k_1} d B_t^{k_2} , \quad  1 \le i \le N,  \\[0.3cm]
& d \langle \mil{x, i}{\cdot}, \mil{x, j}{\cdot} \rangle_t 
=  a^{ij} (x) dt 
+ \tfrac{1}{2} \sum_{\substack{1 \le k_1 \le d \\ 0 \le k_2 \le d}} 
\sigma_{k_1}^i(x) 
\{L_{k_1} \sigma_{k_2}^j (x)  +  L_{k_2} \sigma_{k_1}^j (x)\} B_t^{k_2} dt \\
& \qquad \qquad + \tfrac{1}{2} \sum_{\substack{1 \le k_1 \le d \\ 0 \le k_2 \le d}} 
\sigma_{k_1}^j (z)  \{L_{k_1} \sigma_{k_2}^i (x) +  L_{k_2} \sigma_{k_1}^i (x)\} B_t^{k_2} dt  \\ 
& \qquad \qquad  + \tfrac{1}{4} \sum_{\substack{0 \le k_1, k_2 \le d \\ 1 \le k_3 \le d}}
\bigl\{ L_{k_1} \sigma_{k_3}^i (x) + L_{k_3} \sigma_{k_1}^i (x) \bigr\} 
\bigl\{ L_{k_2} \sigma_{k_3}^j (x) + L_{k_3} \sigma_{k_2}^j (x) \bigr\} B_t^{k_1} B_t^{k_2}  dt, \quad 1 \le i, j \le N. 
\end{align*} 
For the term $T_1$, the stochastic Taylor expansion gives: 
\begin{align*}
T_1 
& = h \cdot \partial_t u (kh, z) 
+ \E \Bigl[ \int_0^h \int_0^s  \partial^2_t u (kh + v, \bar{X}_v^z) 
dv ds \Bigr] 
+ \sum_{i = 1}^N \E \Bigl[ \int_0^h \int_0^s \partial_i \partial_t u (kh + v, \bar{X}^{z}_v) d \bar{X}_v^{z, i} ds \Bigr]  \\ 
& \qquad + \tfrac{1}{2} \sum_{i, j = 1}^N \E \Bigl[ \int_0^h \int_0^s \partial_{ij} \partial_t u (kh + v, \bar{X}^{z}_v) d \langle \bar{X}^{z, i}, \bar{X}^{z, j} \rangle_v  ds \Bigr] \\ 
& = h \cdot \partial_t u (kh, z) +  \tfrac{h^2}{2} \cdot \partial_t^2 u (kh, z) + \tfrac{h^2}{2} \cdot L_0 \partial_t u (kh, z) 
+ \mathscr{R}_{k, 1}^{(n)} (h, z), 
\end{align*}
where $\mathscr{R}_{k, 1}^{(n)} (h, z)$ is the remainder term such that $\mathbb{E} [\mathscr{R}_{k, 1}^{(n)} (h, z)]$ is given in the form of (\ref{eq:err}). 

For the second term $T_2$, we have: 
\begin{align*} 
T_2 
& = \sum_{i = 1}^N \E \Bigl[ \int_0^h \partial_i u (kh + s, \bar{X}_s^z) b^i(z) ds \Bigr]
+ \tfrac{1}{2} \sum_{i = 1}^N \sum_{k = 0}^d \E \Bigl[ \int_0^h \partial_i u (kh + s, \bar{X}_s^z) \bigl\{ 
L_k \sigma_0^i (z) + L_0 \sigma_k^i (z) \bigr\} B_{s}^k ds \Bigr] \\ 
& =  \sum_{i = 1}^N 
\Bigl\{h \cdot \partial_i u (kh, z) +  \tfrac{h^2}{2} \cdot \partial_t \partial_i u (kh, z) + \tfrac{h^2}{2} \cdot L_0 \partial_i u (kh, z)  \Bigr\} b^i (z) + \tfrac{h^2}{2} \sum_{i = 1}^N \partial_i u (kh, z) L_0 \sigma_0^i (z) \\ 
& \quad + \tfrac{h^2}{4} \sum_{i,j = 1}^N \sum_{k_1 = 1}^d \partial_{ij} u (kh, z) \sigma_{k_1}^i (z) 
\bigl\{ L_{k_1} \sigma_0^j (z) +  L_0 \sigma_{k_1}^j (z)  \bigr\} + \mathscr{R}_{k, 2}^{(n)} (h, z),  
\end{align*}
where $\mathscr{R}_{k, 2}^{(n)} (h, z)$ is the remainder term such that $\mathbb{E} [\mathscr{R}_{k, 2}^{(n)} (h, z)]$ is given in the form of (\ref{eq:err}).

Finally, for the term $T_3$ we have: 
\begin{align*} 
T_3 & = \tfrac{1}{2} \sum_{i, j = 1}^N \E 
\Bigl[ \int_0^h \partial_{ij} u (kh + s, \bar{X}_s^z) a^{ij} (z) ds \Bigr] \\
& \quad + \tfrac{1}{2} \sum_{i, j = 1}^N 
\sum_{\substack{1 \le k_1 \le d \\ 0 \le k_2 \le d }} \E 
\Bigl[ \int_0^h \partial_{ij} u (kh + s, \bar{X}_s^z) \sigma_{k_1}^i (z) 
\bigl\{L_{k_1} \sigma_{k_2}^j (z) 
+ L_{k_2} \sigma_{k_1}^j (z)   \bigr\} B_s^{k_2} ds \Bigr]  \\ 
& \quad + \tfrac{1}{8} \sum_{i, j = 1}^N 
\sum_{\substack{0 \le k_1, k_2 \le d \\ 1 \le k_3 \le d }} 
\E \Bigl[ \int_0^h \partial_{ij} u (kh + s, \bar{X}_s^z) \bigl\{L_{k_1} \sigma_{k_3}^i (z) 
+ L_{k_3} \sigma_{k_1}^i (z)   \bigr\}
\bigl\{L_{k_2} \sigma_{k_3}^j (z) 
+ L_{k_3} \sigma_{k_2}^j (z)   \bigr\} B_s^{k_1} B_s^{k_2} ds \Bigr]  \\
& = \tfrac{h}{2}  \sum_{i,j = 1}^N \partial_{ij} u (kh, z) a^{ij} (z) 
+ \tfrac{h^2}{4} \sum_{i, j= 1}^N  \partial_t \partial_{ij} u (kh, z) a^{ij} (z)
+ \tfrac{h^2}{4} \sum_{i, j, l= 1}^N  
\partial_{ijl} u (kh, z) a^{ij} (z) b^l (z) \\
& \quad + \tfrac{h^2}{8} \sum_{i, j, m, l= 1}^N   \partial_{ijml} u (kh, z) a^{ij} (z) a^{ml} (z)
+ \tfrac{h^2}{4} \sum_{i, j = 1}^N 
\sum_{\substack{1 \le k_1 \le d }} \partial_{ij} u (kh, z) \sigma_{k_1}^i (z) 
\bigl\{L_{k_1} \sigma_{0}^j (z) 
+ L_{0} \sigma_{k_1}^j (z)   \bigr\}  \\ 
& \quad + \tfrac{h^2}{2} \sum_{i, j, l = 1}^N 
\sum_{\substack{1 \le k_1 , k_2 \le d }} 
\partial_{ijl} u (kh , z) \sigma_{k_1}^i (z) 
\sigma_{k_2}^l (z) L_{k_1} \sigma_{k_2}^j (z)  \\
& \quad +  \tfrac{h^2}{8} \sum_{i, j = 1}^N \sum_{k_1, k_2 = 1}^d  \partial_{ij} u (kh, z) \bigl\{ 
L_{k_1} \sigma_{k_2}^i (x) L_{k_1} \sigma_{k_2}^j (x) 
+ L_{k_1} \sigma_{k_2}^i (x) L_{k_2} \sigma_{k_1}^j (x) \bigr\} 
+  \mathscr{R}_{k, 3}^{(n)} (h, z),
\end{align*}
where $\mathscr{R}_{k, 3}^{(n)} (h, z)$ is such that $\mathbb{E} [\mathscr{R}_{k, 3}^{(n)} (h, z)]$ is given in the form of (\ref{eq:err}). Here, we introduce the differential operator: for $\varphi \in C_b^\infty (\mathbb{R}^N; \mathbb{R}))$,  $\xi, x \in \mathbb{R}^N$, 
\begin{align}
L_0^\xi  \varphi (x) = \sum_{i = 1}^N b^i (\xi) \partial_i \varphi (x)
+ \tfrac{1}{2} \sum_{i,j = 1}^N a^{ij} (\xi) \partial_{ij} \varphi (x). 
\end{align}
Making use of $\partial_t u = L_0 u$ and the above operator $L_0^\xi$, we obtain: 
\begin{align*}
F_k^{(n)}(z)
& = h \cdot \partial_t u (kh, z) 
+ h \cdot L_0 u (kh, z)  
+ \tfrac{h^2}{2} L_0 \partial_t u (kh, z) 
+ \tfrac{h^2}{2} L_0 L_0^\xi u (kh, z)|_{\xi = z} 
+ \tfrac{h^2}{2} \sum_{i = 1}^N \partial_i u (kh, z) L_0 \sigma_0^i (z)\\
& \quad + \tfrac{h^2}{2} \sum_{i,j = 1}^N \sum_{k_1 = 1}^d \partial_{ij} u (kh, z) \sigma_{k_1}^i (z) \bigl\{ L_{k_1}\sigma_0^j (z) + L_0 \sigma_{k_1}^j (z) \bigr\}  + \Phi_2 (kh, z) + \mathscr{R}_k^{(n)} (h ,z) \\
& = - \tfrac{h^2}{2} L_0^2 u (kh, z) 
+ \tfrac{h^2}{2} L_0 L_0^\xi u (kh, z)|_{\xi = z} 
+ \Phi_1 (kh, z) +  \Phi_2 (kh, z) + \mathscr{R}_k^{(n)} (h ,z),  
\end{align*}
where we have defined
$ \textstyle
\mathscr{R}_k^{(n)} (h ,z) = \sum_{i = 1}^3 \mathscr{R}_{k, i}^{(n)} (h ,z)
$. Since it follows from (\ref{eq:em_err}) that 
\begin{align}
\Phi^{\mrm{EM}} (kh, z) = \tfrac{h^2}{2} L_0^2 u (kh, z) 
- \tfrac{h^2}{2} L_0 L_0^\xi u (kh, z) |_{\xi = z},  
\end{align} 
we have 
\begin{align} \label{eq:final_F}
F_k^{(n)} (z) = - \Phi^{\mrm{EM}} (kh, z) + \Phi_1 (kh, z) 
+ \Phi_2 (kh, z)  + \mathscr{R}_k^{(n)} (h , z)
 = - \Phi_3 (kh, z) + \mathscr{R}_k^{(n)} (h , z).  
\end{align}
%
%
From (\ref{eq:local_err}) and (\ref{eq:final_F}), we conclude. 

\section{Proof of Proposition \ref{prop:small_diff_bd}} \label{appendix:small_diff_bd} 

Before showing Proposition \ref{prop:small_diff_bd}, we prepare some tools from Malliavin calculus. The details are found in, e.g. \cite{nua06}. 

Let $(\Omega, \mathcal{F}, \mu)$ be $d$-dimensional Wiener space and $H = L^2 ([0,T]; \mathbb{R}^d)$ equipped with the inner product $\langle \cdot, \cdot \rangle_H$ defined as $\textstyle \langle h_1, h_2 \rangle_H = \int_0^T h_1 (s) \cdot h_2 (s) ds, \, h_1, h_2 \in H$. For $h \in H$, we write $\textstyle W(h) \equiv \int_0^T h(s) \cdot d W_s$ and introduce a space of some Wiener functionals as 
$\mathscr{S} \equiv \bigl\{ F = f \bigl(W(h_1), \ldots, W(h_n) \bigr) \, | \,  f \in C_p^\infty (\mathbb{R}^n; \mathbb{R}), \, h_1, \ldots, h_n \in H, \, n \ge 1  \bigr\}$. For $p \ge 1$, the Malliavin derivative operator $D : \mathscr{S} \to L^p (\Omega; H)$ is defined via the following operation: 
$$ \mathscr{S} \ni F \mapsto \ 
D_t F = \sum_{i = 1}^n \partial_i f (W(h_1), \ldots, W(h_n)) h_i (t), \ t \in [0, T].
$$
For $F \in \mathscr{S}$, $\{D_t F \}_{t \ge 0}$ is typically treated as a $d$-dimensional stochastic process, and  we write $D_{j ,t} F, \, 1 \le j \le d$ as the $j$-th element of $D_t F$. For $k \in \mathbb{N}$, the $k$-th order Malliavin derivative is denoted by $D^k$ which is closable from $\mathscr{S}$ into $ L^p (\Omega; H^{\otimes k})$. We write $\mathbb{D}_{k, p}$, $k \in \mathbb{N}$, $p \ge 1$ as the completion of $\mathscr{S}$ w.r.t. the norm 
$$
\| F \|_{k, p} 
= \bigl\{ \mathbb{E} [|F|^p] + \sum_{j = 1}^k \mathbb{E} 
[ \| D^j F \|_{H^{\otimes k}}^p] \bigr\}^{1/p}, \ \ F \in \mathscr{S}.  
$$
Define $\textstyle \mathbb{D}^\infty = \bigcap_{k \in \mathbb{N}} \bigcap_{p \ge 1} \mathbb{D}_{k, p}$. For a $m$-dimensional random vector $F = (F^1, \ldots, F^m) \in (\mathbb{D}^\infty)^m$, the Malliavin covariance $\sigma^F = (\sigma_{ij}^F)_{1 \le i,j \le m}$ is defined as $\sigma_{ij}^F = \langle D F^i, D F^j \rangle_H$. Then, $F$ is said to be \emph{non-degenerate in Malliavin sense} if $\sigma^F$ is invertible a.s. and also $\textstyle (\det \sigma^F)^{-1} \in \bigcap_{p \ge 1} L^p (\Omega)$. For non-degenerate Wiener functionals, we have the integration by parts on Wiener space as follows: let $G \in \mathbb{D}^\infty$, $
\varphi \in C_b^\infty (\mathbb{R}^m)$ and $F \in (\mathbb{D}^\infty)^m$ be non-degenerate in Malliavin sense. Then, for any multi-index $\alpha \in \{1, \ldots, m\}^k, \, k \in \mathbb{N}$, there exists $\mathscr{H}_\alpha (F, G)$ such that 
\begin{align} \label{eq:ibp}
\mathbb{E} [\partial_\alpha \varphi (F) G ] 
= 
\mathbb{E} [\varphi (F) \mathscr{H}_\alpha (F, G) ] 
\end{align} 
where the stochastic weight $\mathscr{H}_\alpha (F, G)$ is recursively defined via: 
$$ 
\mathscr{H}_{(\alpha_1)} (F, G)
= \sum_{j = 1}^m \delta \bigl( (\sigma^F)^{-1}_{\alpha_1 j} G \, DF \bigr), 
\quad 
\mathscr{H}_{(\alpha_1, \ldots, \alpha_k)} (F, G)
= \mathscr{H}_{(\alpha_k)} 
\bigl(F, \mathscr{H}_{(\alpha_1, \ldots, \alpha_{k-1})} (F, G) \bigr)
$$
with $\delta$ denoting the adjoint operator of $D$. 
\\

\noindent 
(\textit{Proof of Proposition \ref{prop:small_diff_bd}}.) 
We recall $u^\varepsilon (s, x)  = \E \bigl[  f (X_{T-s}^{\varepsilon,x}) \bigr]$ and set $\alpha \in \{1, \ldots, N \}^k, \, k \in \mathbb{N}$ throughout the proof. We have that
\begin{align}
	\Bigl| \int_0^T \E \bigl[ g (X_s^{x, \varepsilon}) 
	\partial_{\alpha} u^\varepsilon (s, X_s^{x, \varepsilon}) \bigr]  ds 
	\Bigr| = J_1 + J_2, 
\end{align} 
where 
\begin{align}
	J_1 = \Bigl| \int_0^{T/2}  \E \bigl[ g (X_s^{x, \varepsilon}) 
	\partial_{\alpha} u^\varepsilon (s, X_s^{x, \varepsilon}) \bigr]  ds \Bigr|, 
	\qquad  
	J_2 =  \Bigl| \int_{T/2}^T  \E \bigl[ g (X_s^{x, \varepsilon}) 
	\partial_{\alpha} u^\varepsilon (s, X_s^{x, \varepsilon}) \bigr]  ds \Bigr|.  
\end{align} 
We then derive upper bounds for the terms $J_1$ and $J_2$. 
We first consider the term $J_1$. We note that $X_{T-s}^{\xi, \varepsilon}$ is non-degenerate in Malliavin sense under the condition \eqref{con:coeff} and \eqref{con:hor}. Thus, applying the Malliavin integration by parts (\ref{eq:ibp}), we have that: for any $(s, \xi) \in [0, T] \times \mathbb{R}^N$ and $1 \le i \le N$, 
\begin{align}
\partial_{\xi_i} u^\varepsilon (s, \xi ) 
 = \sum_{j = 1}^N \E \bigl[ \partial_j f (X_{T-s}^{\xi, \varepsilon})  J_{T-s, ij}^{\xi, \varepsilon} \bigr] 
 =  \sum_{j = 1}^N \E \bigl[ f (X_{T-s}^{\xi, \varepsilon})
 \mathscr{H}_{(j)}  \bigl( X_{T-s}^{\xi, \varepsilon},  J_{T-s, ij}^{\xi, \varepsilon} \bigr) \bigr],
\end{align}
where $J_{T-s}^{\xi, \varepsilon} = \bigr( J_{T-s, ij}^{\xi, \varepsilon} \bigr)_{1 \le i, j \le N}$ is defined as $J_{T-s, ij}^{\xi, \varepsilon}  = \tfrac{\partial X_{T-s}^{\xi , \varepsilon,  j}}{\partial \xi_i} $. We have that $J_{T-s, ij}^{\xi, \varepsilon} \in \mathbb{D}^\infty$ and $\| J_{T-s, ij}^{\xi, \varepsilon} \|_{k, p} < \infty$ for any $k \in \mathbb{N}$ and $p \ge 1$ due to the condition \eqref{con:coeff}. Then, it holds that  
\begin{align}
\bigl| \partial_\alpha u^\varepsilon (s, \xi )  \bigr| 
& \le 
\tfrac{K(T)}{(T-s)^q \varepsilon^{k}} \| f \|_\infty  ( 1 + |\xi|^Q)  
\label{eq:bd_grad}
\end{align} 
for some constants $q, Q  > 0, \, p \ge 1 $ and non-decreasing function $K(\cdot)$ independent of $x$ and $\varepsilon$. To obtain the bound (\ref{eq:bd_grad}), we have exploited the following estimate whose proof is postponed to the end of this section: 
\begin{lemma} \label{lemma:bd_weight}
Let $x \in \mathbb{R}^N$ and $G \in \mathbb{D}^\infty$ satisfying that: for any $j \in \mathbb{N}$ and $q \ge 1$, there exists a constant such that  $\| G \|_{j, p} \le c$ for some $c > 0$. Assume the conditions \eqref{con:coeff} and \eqref{con:hor} with $M=1$ hold.  Then, for any $\beta \in \{1, \ldots, N\}^k, \, k \in \mathbb{N}$ and $p \ge 1$,  there exist constants $q, Q > 0$ and a non-decreasing function $K(\cdot)$ independent of $x$ and  $\varepsilon \in (0,1)$ such that 
\begin{align}
\bigl\{ \E \bigl[  \bigl| 
\mathscr{H}_{(\beta)} (X_{t}^{x, \varepsilon}, G ) \bigr|^p  \bigr] \bigr\}^{1/p} 
\le  \frac{K(t)}{t^q \varepsilon^k} (1 + |x|^Q). 
\end{align}
\end{lemma}
Using the bound (\ref{eq:bd_grad}) and the polynomial growth of $g$ with moment bound for $X_s^x$ under (\ref{con:coeff}), we have that: 
\begin{align}
J_1 
\le \frac{K (T) }{T^{q} \varepsilon^{k}}   \| f \|_\infty  \times \sup_{s \in [0, T]} \mathbb{E} [ (1 + |X_{T-s}^{x, \varepsilon} |^{Q_1}) ] 
\le \frac{K (T) }{T^{q} \varepsilon^{k}}   \| f \|_\infty  (1 + |x|^{Q_2}) 
\end{align}
for some constants $C_1, q, Q_1, Q_2 > 0$ and non-decreasing functions $K (\cdot)$ which are independent of $x$ and $\varepsilon$. For the second term $J_2$, we apply the Malliavin integration by parts and Lemma \ref{lemma:bd_weight} to obtain that:
\begin{align}
J_2 
& =  \Bigl| \int_{T/2}^T  \E \bigl[ 
u^\varepsilon (s, X_s^{x, \varepsilon})
\mathscr{H}_{(\alpha)} \bigl(  X_s^{x, \varepsilon}, g (X_s^{x, \varepsilon}) \bigr) \bigr]  ds \Bigr|  \nonumber \\ 
& \le   \| f \|_\infty \int_{T/2}^T 
\E \bigl[  \bigl|  \mathscr{H}_{(\alpha)} \bigl(  X_s^{x, \varepsilon}, g (X_s^{x, \varepsilon}) \bigr)   \bigr| \bigr] \, ds 
\le \frac{K(T)}{T^q  \varepsilon^{k}}
\| f \|_\infty  ( 1 + |x|^Q), 
\end{align}  
for some constants $q, Q > 0$ and a non-decreasing function $K$ independent of $x \in \mathbb{R}^N$ and $\varepsilon \in (0,1)$, where in the last line we have used the polynomial growth of $g$ and the following estimate under the condition \eqref{con:coeff}: for any $k \in \mathbb{N}$ and $p \ge 1$ there exist constants $c, Q > 0$ independent of $s, x, \varepsilon$ such that $\| g (X_s^{x, \varepsilon}) \|_{k, p} < c (1 + |x|^Q)$. The proof of Proposition \ref{prop:small_diff_bd} is now complete.   
\subsection{Proof of Lemma \ref{lemma:bd_weight}}
%
We have from (2.32) in \citep{nua06}  that: for $1 \le p \le q < \infty$, 
\begin{align}
\bigl\{ \E \bigl[ |  \mathscr{H}_{(\beta)} (X_t^{x, \varepsilon},  G) |^p \bigr]
\bigr\}^{1/p} 
\le C_1 \bigl\|  \gamma^{X_t^{x, \varepsilon}}  D X_t^{x, \varepsilon} \bigr\|_{k, \gamma}^k  \| G \|_{k, q} 
\le C_2 \bigl\|  \gamma^{X_t^{x, \varepsilon}}  D X_t^{x, \varepsilon} \bigr\|_{k, \gamma}^k, 
\end{align}
with some constants $C_1, C_2 > 0$ and $\gamma > 1$. The Malliavin derivative of $X_t^{x, \varepsilon, i}, \, 1 \le i \le N$,  is given as (for instance see \citep{nua06, shige04}):
\begin{align} \label{eq:MD}
D_{j, s} X_t^{x, \varepsilon, i} =
\varepsilon 
\bigl[  J_t^{x, \varepsilon} \bigl( J_s^{x, \varepsilon} \bigr)^{-1} \sigma_j (X_s^{x, \varepsilon}) \bigr]_i,  \qquad s \in [0, t], \ \  1 \le j \le d.   
\end{align}
We define $N\times N$ matrix valued random variable $\widetilde{M}^{X_t^{x, \varepsilon}} = \bigl( \widetilde{M}^{X_t^{x, \varepsilon}}_{i_1 i_2} \bigr)_{1 \le i_1, i_2 \le N} $ as: 
\begin{align}
\widetilde{M}^{X_t^{x, \varepsilon}} 
= \sum_{j = 1}^d \int_0^t  J_t^{x, \varepsilon} 
\bigl( J_s^{x, \varepsilon} \bigr)^{-1} V_j (X_s^{x, \varepsilon}) \,  \otimes J_t^{x, \varepsilon} \bigl( J_s^{x, \varepsilon} \bigr)^{-1} V_j (X_s^{x, \varepsilon}) ds. 
\end{align} 
Under the uniform H\"ormander condition \eqref{con:hor} with $M=1$, the matrix $\widetilde{M}^{X_t^{x, \varepsilon}}$ is shown to be invertible a.s. and also $\mathbb{E} \bigl[ \bigl| (\det \widetilde{M}^{X_t^{x, \varepsilon}})^{-1} \bigr|^p \bigr] < \infty$ for all $p \ge 1$. 
We then write the inverse of the matrix $\widetilde{M}^{X_t^{x, \varepsilon}}$ as $\widetilde{\gamma}^{X_t^{x, \varepsilon}}$.  
By considering stochastic Taylor expansion of $\bigl( J_s^{x, \varepsilon} \bigr)^{-1} V_j (X_s^{x, \varepsilon})$, it is shown from  Theorem 6.16 (and its proof) in \citep{shige04} that: there exist constants $q, Q > 0$ and a non-decreasing function $K(\cdot)$  independent of $x, \varepsilon, t$ such that 
\begin{align}
\bigl\| \widetilde{\gamma}^{X_t^{x, \varepsilon}}_{i_1 i_2} \bigr\|_{m, n}  \le \frac{K(t)}{t^q} (1 + |x|^Q), \qquad m, n \in \mathbb{N}, \quad 1 \le i_1, i_2 \le N.  
\end{align}
In particular, noticing that the Malliavin covariance of $X_{t}^{x, \varepsilon}$ is given as $M^{X_t^{x, \varepsilon}} = \varepsilon^2 \widetilde{M}^{X_t^{x, \varepsilon}}$, we have that: for $m, n \in \mathbb{N}$ and $1 \le i_1, i_2 \le N$, 
\begin{align} \label{eq:bd_inv}
\bigl\|  \gamma^{X_t^{x, \varepsilon}}_{i_1 i_2} \bigr\|_{m, n} 
=  \varepsilon^{- 2}
\bigl\| \widetilde{\gamma}^{X_t^{x, \varepsilon}}_{i_1 i_2}  \bigr\|_{m, n} 
\le \frac{K(t)}{t^q \varepsilon^{2}} (1 + |x|^Q). 
\end{align}
We thus obtain from (\ref{eq:MD}), (\ref{eq:bd_inv}) and H\"older's inequality that:  
\begin{align}
\bigl\|  \gamma^{X_t^{x, \varepsilon}}  D X_t^{x, \varepsilon} \bigr\|_{k, \gamma}^k 
\le \frac{K(t)}{t^q \varepsilon^{k}} (1 + |x|^Q) ,
\end{align} 
and we conclude. 
\bibliographystyle{rss}
\bibliography{milstein} 

\end{document}